\documentclass[12pt]{article}

\usepackage{amsmath}
\usepackage{amsthm}
\usepackage{mathtools}
\usepackage{dashbox}
\usepackage{a4wide}
\usepackage{epsfig}
\usepackage{theorem}
\usepackage{amssymb}
\usepackage{bbm}
\usepackage[T1]{fontenc}
\usepackage{cancel}
\usepackage{blkarray}
\usepackage{multirow}
\usepackage{tikz}
\usepackage{makeidx}
\usepackage{graphicx}
\usepackage{subfig}

\newcommand{\cL}{{\cal L}}

\newcommand{\cF}{{\cal F}}

\newcommand{\bi}{\bigskip}
\newcommand{\no}{\noindent}
\newcommand{\bea}{\begin{eqnarray}}
\newcommand{\eea}{\end{eqnarray}}

\newcommand{\be}{\begin{equation}}
\newcommand{\ee}{\end{equation}}

\newcommand{\lk}{\left(}\usepackage{multirow}

\newcommand{\sli}{\sum\limits}
\newcommand{\pli}{\prod\limits}

\newcommand{\vk}{\vec{k}}

\newcommand{\il}{\int\limits}

\newcommand{\ZZ}{\mathbbm{Z}}
\newcommand{\RR}{\mathbbm{R}}
\newcommand{\NN}{\mathbbm{N}}
\usepackage{amsmath}
\usepackage{a4wide}
\usepackage{epsfig}
\usepackage{theorem}
\usepackage{amssymb}
\usepackage{bbm}
\usepackage[latin1]{inputenc}

\allowdisplaybreaks[3]

\begin{document}

\title{On Riemann's Paper, ``On the Number of Primes Less Than a Given Magnitude''}

\author{W.~Dittrich\\
Institut f\"ur Theoretische Physik\\
Universit\"at T\"ubingen\\
Auf der Morgenstelle 14\\
D-72076 T\"ubingen\\
Germany\\
electronic address: qed.dittrich@uni-tuebingen.de
}
\date{\today}

\maketitle
\bi

\no

\begin{abstract}
This paper is devoted to one of the members of the G\"ottingen triumvirate, Gau\ss{}, Dirichlet and Riemann. 
It is the latter to whom I wish to pay tribute, and especially to his world-famous article of 1859, which he presented in person 
at the Berlin Academy upon his election as a corresponding member. His article entitled, 
``\"Uber die Anzahl der Primzahlen unter einer gegebenen Gr\"o\ss{}e'' (``On the Number of Primes Less Than a Given Magnitude''),
revolutionized mathematics worldwide. Included in the present paper is a detailed analysis of Riemann's article, 
including such novel concepts as analytical continuation in the complex plane; 
the product formula for entire functions; and, last but not least, a detailed study of the zeros of the so-called
Riemann zeta function and its close relation to determining the number of primes up to a given magnitude, 
i.e., an explicit formula for the prime
counting function. 
\end{abstract}

\section*{Short Biography of Bernhard Riemann (1826 - 1866)}

Bernhard Riemann was born in Breselenz near Dannenberg in Lower Saxony in 1826. Like his father, he was first supposed
to become a pastor, but already in high school Riemann's extraordinary mathematical talent caught the attention of his principal.
It is said that Riemann read the 859-page book by Legendre on number theory which was loaned to him by the principal in one week.
He began studying mathematics in G\"ottingen, where he attended lectures by Gau\ss{}, although they were closed to first-semester 
students. Riemann then transferred to Jacobi and Dirichlet in Berlin,
both of whom supported and encouraged him;
he then returned to G\"ottingen. His doctoral thesis was on the theory of functions.
In order to be permitted to teach as private lecturer in G\"ottingen,
candidates had to submit three suggestions for the topic of their Habilitation lecture,
and normally the department head would choose the uppermost topic on the list.
Riemann's third topic was ``Basics of Geometry'', and when Gau\ss{} read that, he, as department head,
selected that topic for Riemann's Habilitation lecture. Very much surprised, 
Riemann put all his research on the topic ``Electricity, Magnetism, Light and Graviation'' aside and in 1854, two months before his 
trial lecture, created the foundations of differential geometry.
Gau\ss{} was thrilled! In 1855 Gau\ss{} died and was followed by Dirichlet. 
When Dirichlet died four years later, Riemann took over the mathematics chair
at G\"ottingen University. In 1862, he married Elise Koch, with whom he had one daughter. 
Riemann fell ill with TB and looked for relief in the milder climate 
of Tessin, where he died at the early age of only 39 at Lake Maggiore.

In addition to founding differential geometry, Riemann made other major contributions; 
especially important was his work in the theory of functions;
his ``\"Uber die Anzahl der Primzahlen unter einer gegebenen Gr\"o\ss{}e'' (On the Number of Primes up to a Given Magnitude), communicated
in the ``Monatsberichte der Berliner Akademie, November 1859, with findings on the zeta function; his works on the theory of 
integration, the Fourier transformation, the hypergeometric differential equation,
and the hyperbolic differential equations and stability problems of solutions of partial differential equations in mathematical physics.
Riemann was influenced by the research on algebraic geometry and topology by his Italian mathematician friends Betti and Beltrami. 
Einstein's General theory of Relativity would be unthinkable without Riemannian Geometry. 

All these topics have kept mathematicians and theoretical physicists busy for many years 
and will continue to do so for many more to come. Today, exactly 150 years after Riemann's 
death, the major unsolved problem in pure mathematics is the so-called Riemann hypothesis,
a conjecture made by Riemann in 1859 in his paper on the number of primes less than a 
given positive integer $x$. 

Mathematicians later realized that Riemann's hypothesis governs the distribution of prime 
numbers to an extraordinary extent, which is why its proof is so eagerly sought. 
Since all the efforts of some of the best mathematicians have failed so far, 
perhaps another Riemann is needed. 

This is also true for many local relativistic quantum field theory models of elementary 
particles, where Riemann's results are of utmost importance for handling infinities 
with the aid of his zeta-function regularization. In non-relativistic quantum mechanics, 
we need a Riemannian Hamiltonian which becomes diagonalized in the prime number basis. 
The measurement process, i.e., the operator acting on an object that will provide 
us with such a set of discrete prime number eigenvalues, is still to be found. 
One might wonder what kind of symmetry structure lies behind this kind of physical system.

Let us not forget that the few papers that Riemann published in his lifetime dealt with 
physics problems. Moreover, in the days of Gauss, Dirichlet and Riemann, a distinction 
between the disciplines of mathematics and physics was not made. In particular, 
Riemann approached problems in mathematics and physics not so much as an analyst but  
illuminated them globally from a geometric and topological viewpoint, meaning that he 
made many results of analysis better understood using the new methods of the theory of 
functions and analytical continuation into the whole complex plane, thereby simplifying 
many problems of real analysis.

\section{Towards Euler's Product Formula and Riemann's Extension of the Zeta Function}

There is a very close connection between the sums of the reciprocals of the integers raised to a variable power that 
Euler wrote down in 1737, the now-called zeta function,
\be
\zeta (s) = \sli^\infty_{n = 1} \frac{1}{n^s} = 1 + \frac{1}{2^s} + \frac{1}{3^s} + \frac{1}{4^s} + \frac{1}{5^s} + \cdots \, , 
\quad \quad s > 1
\ee
and the primes - which, as integers, are the very signature of discontinuity. Euler considered $s$
to be a real integer variable with $s >1$ to insure convergence of the sum. Multiplying the definition of $\zeta (s)$ by ${1/2}^s$ we obtain
\be
\frac{1}{2^s} \zeta (s) = \sli^\infty_{n = 1}  \frac{1}{(2n)^s} = \frac{1}{2^s} + \frac{1}{4^s} +  \frac{1}{6^s} + \frac{1}{8^s} + \cdots \, 
\ee
and subtracting this from $\zeta (s)$ we get
\begin{align}
 \zeta (s) - \frac{1}{2^s} \zeta (s) &=  \sli^\infty_{n = 1}  \frac{1}{n^s}  - \sli^\infty_{n = 1}  \frac{1}{(2 n)^s} \nonumber\\
 \mbox{or} \quad \quad  \lk 1 - \frac{1}{2^s} \right) \zeta (s) &= 1 + \frac{1}{3^s} + \frac{1}{5^s} + \frac{1}{7^s} + \frac{1}{9^s}  +
 \frac{1}{11^s} + \cdots \, .
 \end{align}
 Hence all the multiples of the prime $n = 2$ disappeared from the original sum of the defined $\zeta (s)$. In short, we found
 \be
 \lk 1- \frac{1}{2^s} \right) \zeta (s) = \sli^\infty_{n = 1 \atop \Lambda n \neq 2 k} 
 \frac{1}{n^s} \, .
 \ee
Next, we multiply this last result by $1/3^s$ to obtain
\be
 \frac{1}{3^s} \lk 1 - \frac{1}{2^s} \right) \zeta (s) = 
 \sli^\infty_{n = 1 \atop \Lambda n \neq 2 k}  \frac{1}{(3n)^s} = 
 1 + \frac{1}{3^s} + \frac{1}{9^s} +  \frac{1}{15^s} + \frac{1}{21^s} + \cdots
\ee
and so, subtracting this from $(1- 1/2^s) \zeta (s)$, we have
\begin{align}
 \lk 1 - \frac{1}{2^s} \right)  \lk 1 - \frac{1}{3^s} \right)  \zeta (s) &= 1 + \frac{1}{5^s} + \frac{1}{7^s} +  \frac{1}{11^s} + \cdots \nonumber\\
 &= \sli^\infty_{
 \begin{array}{c}
                  n = 1 \\
                  \Lambda n \neq 2 k \\
                  \Lambda n \neq 3 k
                 \end{array}
} \frac{1}{n^s} \, .
\end{align}
Now we multiply this result by $1/5^s$
and so on. As we  repeat this process over and over, multiplying through our last result by $1/p^s$, 
where $p$ denotes successive primes, we subtract out all the multiples of the primes. Hence, since all integers are composed of 
primes (Euclid's fundamental theorem of the theory of numbers), we removed all numbers of the right-hand side of the defining sum of 
$\zeta (s)$ - except for the number 1. Thus our final result is the product
\be
\left\{ \Pi_{p \, \mbox{\scriptsize prime}} \lk 1  - p^{- s} \right) \right\} \zeta (s) = 1
\ee
or
\be
\boxed{\zeta (s) = \Pi_{p \, \mbox{\scriptsize prime}} \frac{1}{1 - p^{- s}} = \sli^\infty_{n = 1} 
\frac{1}{n^s} \, , \quad \quad s > 1 } \, .
\ee

Euler's actual statement reads: ``Si ex serie numerorum primorum sequens formetur expressio $\pli_p \frac{p^s}{(p^s - 1)}$ erit 
eius valor aequalis summae huius seriei $\sli_{n = 1} \frac{1}{n^s}$.''

Now we are going to extend Euler's zeta function into the complex plane $C$, 
which is a major achievement of Riemann's. Hence from now on, s is complex valued and we write
\be
 \zeta (s) = \sli^\infty_{n = 1} \frac{1}{n^s} = \frac{1}{1^2} +  \frac{1}{2^2} +\frac{1}{3^2} + \cdots \quad \quad
\mbox{but with} \quad Re (s) > 1 \, .
\ee
This is an absolutely convergent infinite series, which also holds true for the product of all primes in
\be
\zeta (s) = \Pi_{p \, \mbox{\scriptsize prime}} \frac{1}{1  - p^{- s}} = 
\lk \frac{1}{1 - 2^{-  s}} \right) \cdot 
\lk \frac{1}{1 - 3^{-  s}} \right) \cdot 
\lk \frac{1}{1 - 5^{-  s}} \right) \cdots 
\lk \frac{1}{1 - p^{-  s}} \right) \cdots \, .
\ee
$\zeta (s)$ has no zeros in the region $Re(s)>1$, 
as none of these factors have zeros. However, with Riemann's extension of zeta into the entire complex plane, we will be able to 
locate zeros as well as poles. To show this, we have to analytically continue Euler's original real valued zeta function into 
the entire complex $s$ plane. A first result in this direction will be achieved with the aid of the so-called Dirichlet series, 
which turns up when calculating
\begin{align}
 (1 - 2^{1 - s} ) \zeta (s) &= \sli^\infty_{n = 1} n^{- s} - 2^{1 - s} 
 \sli^\infty_{n = 1} n^{-s} = \sli^\infty_{n = 1} n^{- s} - 2 \sli^\infty_{n = 1}
 (2 n)^{- s} \nonumber\\
 &= 1 - \frac{2}{2^s}   + \frac{1}{2^s} - \frac{2}{4^s} + \cdots = 1 - 
 \frac{1}{2^s} +
  \frac{1}{3^s} - \frac{1}{4^s} + \frac{1}{5^s} - \frac{1}{6^s} + \cdots \nonumber\\
  &= \sli^\infty_{n = 1} \frac{(-1)^{n + 1}}{n^s} =: \eta (s) \, , \quad \quad \mbox{Dirichlet series} \, .
\end{align}
This series is convergent for all $s \in C$  with $ Re(s)>0$.
Hence we can define
\be
\boxed{\zeta (s) = \frac{1}{1 - 2^{1 - s}} \sli^\infty_{n = 1} \frac{(-1)^{n + 1}}{n^s} \quad \quad \mbox{for} \quad Re (s) > 0 \quad 
\mbox{and} \quad 1 - 2^{1 - s} \neq 0 \, .}
\ee
When we write
\begin{align}
 \eta (s) + \frac{2}{2^s} \zeta (s) &= \sli^\infty_{n = 1} \frac{(-1)^{n + 1}}{n^s} + \frac{2}{2^s} \sli^\infty_{n = 1} \frac{1}{n^s} \nonumber\\
 &= \sli^\infty_{n = 1} \lk \frac{1}{(2 n - 1)^s} - \frac{1}{(2 n)^s} + \frac{2}{(2 n)^s} \right) \nonumber\\
 &= \sli^\infty_{n = 1}  \frac{1}{n^s} = \zeta (s) \, ,
\end{align}
we can collect our results so far in the string of formulae
\be
\zeta (s) = \frac{1}{1 - 2^{1 - s}} \sli^\infty_{n = 1} \frac{(-1)^{n + 1}}{n^s} = \frac{\eta (s)}{1 - 2^{1 - s}} = 
\frac{1}{s - 1}  \sli^\infty_{n = 1}  \lk \frac{n}{(n + 1)^s} - \frac{n - s}{n^s} \right) \, .
\ee
Most important, we can continue $\zeta (s)$ into the realm of the critical strip $0<Re(s)<1$.
Of course, the zeros in the denominator in the representation given above have to be excluded, i.e.,  
from
\be
1 - 2^{1 - s} = 0
\ee
follows
\be
1 = e^{(1 - s) \log 2}
\ee
meaning
\be
2 \pi i n = (1 - s) \log 2
\ee
or
\be
s = 1 - \frac{2 \pi i n}{\log 2} \, , \quad \quad n \in \ZZ \, .
\ee
Having shown that the zeta function can be analytically continued into the half plane
$\{ s \in C | Re(s)>0, s \neq 1 \}$, 
we still have to prove that $\zeta (s)$ has a pole at $s=1$:
\begin{align}
 \lim\limits_{s \to 1} \zeta (s) &=  \lim\limits_{s \to 1}  \frac{(s - 1)}{1 - 2^{1 - s}} \sli^\infty_{n = 1}  
 (-1)^{n + 1} n^{- s} = \lim\limits_{s \to 1}  \frac{(s - 1)}{1 - 2^{1 - s}} \log 2 \nonumber\\
  &=  \lim\limits_{s \to 1} \frac{1}{- \log 2 \cdot 2^{1 - s} \cdot (- 1)} \log 2 =  \lim\limits_{s \to 1} \frac{1}{2^{1 - s}} = 1 \, ,
\end{align}
where we used Abel's theorem $\lim_{x \to 1^-} \log(x + 1) = \log 2$ and the continuity of $\log(x + 1)$. How about arguments for the zeta function equal to or less than zero? Later we will show
that the zeta function satisfies the functional equation
\be
\zeta (s) = 2^s \pi^{s - 1} \sin \lk \frac{\pi s}{2} \right) \Gamma (1 - s) \zeta (1 - s) \, .
\ee
This defines $\zeta (s)$ in the whole complex $s$ plane. Note that the left-hand side
goes over by just changing $s \to 1 - s$ into $\zeta (1 - s)$, so that we can compute $\zeta (1 - s)$, 
given $\zeta (s)$, e.g., $\zeta (- 15)$ in terms of $\zeta (16)$.
\bi

\no
If $s$ is a negative even integer, then $\zeta (s) = 0$ because the factor $\sin (\pi s / 2)$ vanishes. These 
are  the trivial zeros of the zeta function. So all non-trivial zeros 
lie in the critical strip where $s$ has a real part between
$0$ and $1$.
\bi

\no
Here is a first curiosity that needs further interpretation.
If one substitutes in the functional equation $s = - 1$, one obtains
\be
\zeta (-  1) = 2^{- 1} \pi^{- 2} (- 1) \Gamma (2) \zeta (2) = \frac{1}{2} \cdot \frac{1}{\pi^2}
( - 1) \cdot 1 \cdot \frac{\pi^2}{6} = - \frac{1}{12} \, ,
\ee
which means that $\zeta (- 1) = - 1/12$.
\bi

\no
This regularized value of $\zeta (- 1)$ has absolutely nothing to do with the real-space representation
of $\zeta (- 1)$ by the divergent series $\zeta (- 1) = \sli^\infty_{n = 1} \frac{1}{n^{- 1}} = 1 + 2 + 3 + 4 + \cdots$, which 
tells us that the same function can have different representations.
Some very learned mathematicians entertain the opinion that the zeta-function regularization has swept away the 
ugly infinities and produced the ``golden nugget'' of the otherwise nonconvergent series. 
In quantum field theory one observes
the same phenomena, where the zeta-function regularization makes infinities disappear (Casimir effect, quantum electrodynamics,
quantum chromodynamics and particle production near black holes). We will come back to this point toward the end of this article. 
\bi

\no
\section{Prime Power Number Counting Function}
\bi

\no
On the way to showing the significance of the zeta zeros for counting prime numbers up to a given magnitude, Riemann introduces an 
important weighted prime number function $f(x)$. We will call it $\Pi (x)$
while others use $J(x)$. Since this function is of utmost importance, we will start 
introducing it by way of examples.

First, the definition of $\Pi (x)$ is given by
\be
\Pi (x) = \sli_{p^n < x \atop p \, \mbox{\scriptsize prime}} \frac{1}{n} \, ,
               \ee
i.e., for every prime number power $p^n$ which is smaller than $x$, we sum up its fractions; for example,
\begin{align}
\Pi (20) = & \lk \frac{1}{1} + \frac{1}{2} + \frac{1}{3} + \frac{1}{4} \right)   
 +
& \lk  \frac{1}{1} + \frac{1}{2} \right) 
 + &
& \lk  \frac{1}{1}\right)  
 + &
& \lk  \frac{1}{1}\right)  
 + &
& \lk  \frac{1}{1}\right) 
 + &
 & \lk  \frac{1}{1}\right)   \nonumber\\
& 
 2^1, 2^2, 2^3, 2^4 < 20 
&
 3^1, 3^2 < 20 
& &
 5^1 < 20 
& &
 7^1 < 20 
& &
 11^1 < 20 
& &
 13^1 < 20 \nonumber\\
& +\lk  \frac{1}{1}\right)  +
\lk  \frac{1}{1}\right)  \nonumber\\
& 17^1 < 20 \quad 19^1 < 20
\end{align}
The brackets can also be reorganized like this:
\begin{align}
 \Pi (20) &= \lk \frac{1}{1} + \frac{1}{1} +\frac{1}{1} +\frac{1}{1} +\frac{1}{1} +\frac{1}{1} +\frac{1}{1} +\frac{1}{1} \right) \nonumber\\
 & + \frac{1}{2} \lk \frac{1}{1} + \frac{1}{1} \right) + \frac{1}{3} \lk \frac{1}{1} \right) + 
  \frac{1}{4} \lk \frac{1}{1}\right) \, .
\end{align}
The first pair of brackets counts the number of primes smaller than $x=20$; 
the second pair counts the primes that are smaller than the square root of $x$, etc. 
Hence, denoting the number of primes up to $x$ by $\Pi (x)$, we get Riemann's formula,
\be
\Pi (x) = \sli^\infty_{n = 1} \frac{1}{n} \pi (x^{1/n}) \, ,
\ee
which contains a finite number of terms, which becomes evident by looking at the following example:
\begin{align}
 \Pi (x) &= \pi (x) + \frac{1}{2} \pi (\sqrt[2]{x}) + \frac{1}{3} \pi  (\sqrt[3]{x}) + 
 \frac{1}{4} \pi  (\sqrt[4]{x})  + \cdots \nonumber\\
 x &= 100 : \nonumber\\
 \sqrt[2]{x} &= 10,  \sqrt[3]{x} = 4.6415 \, ,  
 \sqrt[4]{x} = 3.1622,  \sqrt[5]{x} = 2.51188 \, , \nonumber\\ 
 \sqrt[6]{x} &= 2.15 \ldots,  \sqrt[7]{x} = 1.930 \ldots < 2 \, . 
 \end{align}
If the argument of $\Pi$ is less than $2$, then 
$\Pi (x)=0$. So our result for $\Pi (100)$ is given by
\begin{align}
 \Pi  (100) &= \pi (100) + \frac{1}{2} \pi  (10) + \frac{1}{3} \pi (4.6415) + \frac{1}{4} \pi (3.1622) \nonumber\\
 &  + \frac{1}{5} \pi (2.5118) + \frac{1}{6} \pi (2.15) + 0 + 0 + \cdots \, .
\end{align}
Counting the primes, we obtain
\begin{align}
 \Pi (100) &= 25 + \frac{1}{2} \cdot 4  + \frac{1}{3} \cdot 2 + \frac{1}{4} \cdot 2 + \frac{1}{5} \cdot 1 + \frac{1}{6} \cdot 1 \nonumber\\
 &= 28 \frac{8}{15} = 28.533 \, .
\end{align}
Hence, for any argument $x>1$, the value $\Pi (x)$ can be worked out for a finite sum.
So far we have learned that $\Pi (x)$ measures primes. Evidently $\Pi (x)$ is a step function which starts at 
$\Pi (0) = 0$ and jumps at positive integers, i.e., the jump is $1$ at primes, $1/2$ at squares of primes, and $1/3$ at cubes 
of primes. Hence, our defining equations for $\Pi (x)$ can also be written as
\be
\Pi (x) = \sli_p \sli^\infty_{n = 1} \frac{1}{n} \Theta (x - p^n) \, , 
\ee
where $\Theta (x)$ is the Heaviside step function given by $\Theta (x) = \left\{ \begin{array}{ccl}
                                                                                     1 & , & x > 0 \\
                                                                                     \frac{1}{2} & , & x = 0 \\
                                                                                     0 & , & x < 0
                                                                                    \end{array}
                                                                                    \right. \, .$
\bi

\no
There is still another function of the analytical theory of numbers which we need. It is the so-called 
M\"obius function, which defines the inverse of the zeta function:
\be
\frac{1}{\zeta (s)} = \sli^\infty_{n = 1} \frac{\mu (n)}{n^s} = 1 - \frac{1}{2^s}  - \frac{1}{3^s} - \frac{1}{5^s} +
\frac{1}{6^s} - \frac{1}{7^s}  + \cdots \, .
\ee
Using the original representation
\be
\frac{1}{\zeta (s)} = \lk 1 - \frac{1}{2^s} \right) \lk 1 - \frac{1}{3^s} \right) \lk 1 - \frac{1}{5^s} \right) \lk 1 - \frac{1}{7^s} \right) 
\cdots \, ,
\ee
we may execute the multiplication of the various factors and so end up again with
\be
1 - \frac{1}{2^s} - \frac{1}{3^s} - \frac{1}{5^s} + \frac{1}{6^s} - \frac{1}{7^s} + \frac{1}{10^s} -  \cdots \, ,
\ee
which identifies the following values for $\mu$:
\begin{align}
\mu(1) & = 1, \mu(2) = -1, \mu(3) = -1, \mu(4) = 0, \mu(5) = -1 \, , \nonumber\\
\mu(6) & = 1, \mu(7) = -1, \mu(8) = 0, \mu(9) = 0, \mu(5) = 1, etc.
\end{align}
Here is the rule:
\be
 \mu (n) = \left\{ \begin{array}{ll}
                    - 1 & \mbox{if} \, $n$ \,  \mbox{contains an odd number of primes}\\
                    1 & \mbox{if} \, $n$ \, \mbox{contains an even number of primes}\\
                    0   & \mbox{if} \, $n$ \,  \mbox{contains a quadratic prime factor} 
                   \end{array} \right.
 \ee
For example:
\be
\begin{array}{ll}
\mu (7) = -1; 7 & \mbox{is a prime number}\\
	\mu (66) = -1; 66 = 2 \cdot 3 \cdot 11, & \mbox{odd number of primes}\\
	\mu (18) = 0; 18 = 2 \cdot 3^2, & \mbox{one quadratic prime number}
 \end{array} 
 \ee
For further use we list some lower M\"obius numbers:\\[3em]
\begin{tabular}{|l|l|} \hline
$\mu (n)$ = $- 1$ & $2$ \, $3$ \, $5$ \, $7$ \, $11$ \, $13$ \, $17$ \, $19$ \, $23$ \, $29$ \, $30$ \, $31$ \, $37$ \\ 
\hline
$\mu (n)$ = $0$ & $4$ \, $8$ \, $9$ \, $12$ \, $16$ \, $18$ \, $20$ \, $24$ \, $25$ \, $27$ \, $28$ \, $32$ \, $36$ \\
\hline
$\mu (n) 0 + 1$ & $1$ \, $6$ \, $10$ \, $14$ \, $15$ \, $21$ \, $22$ \, $26$ \, 
$33$ \, $34$ \, $35$ \, $38$ \, $39$\\
\hline
\end{tabular} 
\bigskip

\begin{tabular}{|l|rrrrrrrrrrrrrrrrrrrr|} \hline
$n$ & $1$ &  $2$ & $3$ & $4$ &  $5$ & $6$ & $7$ & $8$ & $9$ & $10$ & $11$ & $12$ & $13$ &
$14$ & $15$ & $16$ & $17$ & $18$ & $19$ &  $20$
\\ 
\hline
$\mu (n)$ & $1$ & $-1$ & $-1$ & $0$ & $-1$ & $1$ & $-1$ & $0$ & $0$ & $1$ & $-1$ & $0$ & $-1$ 
 & $1$ & $1$ & $0$ & $-1$ & $0$ & $-1$ & $0$
\\
\hline
\end{tabular}
\bigskip

 The relation between $\Pi (x)$ and $\pi (x)$ is inverted by Riemann by means of the M\"obius inversion formula to obtain
 \be
 \pi (x) = \sli^\infty_{n = 1} \frac{\mu (n)}{n} \Pi (x^{1/n}) = \Pi (x) - \frac{1}{2} \Pi (x^{1/2}) - \frac{1}{3} 
 \Pi (x^{1/3}) - \frac{1}{5} \Pi (x^{1/5}) + \frac{1}{6} \Pi (x^{1/6}) + \cdots \, .
 \ee
 In the final part of this section I want to discuss briefly a certain integral 
 transform which will be of great help in the next chapter.
 This transformation with kernel $K(z, \xi) = \xi^{z-1}$ 
 is known as Mellin transform, although Riemann knew about it forty years before it 
 became known under this name.
\bi

\no
Let us start with 
 \be
 g  (z) = \il^\infty_{0} d \xi \xi^{z - 1} f (\xi) \, ,
 \ee
for example, with the left-hand side given by $\Gamma (s), Re(s)>0$ and $f(x) = e^{-x}$:
\be
\Gamma (s) = \il^\infty_0 d x e^{- x} x^{s - 1} \, \mbox{with inverse}  \,
e^{- x} = \frac{1}{2 \pi i} \il^{a + i_\infty}_{a - i_{\infty}} d s \frac{\Gamma (s)}{x^s} \, .
\ee
Now we replace $x$ by $nx (n = 1,2,3 ...)$, 
then multiply the equations by constants $c_n$ and sum over $n$:
\begin{align}
 \sli^\infty_{n = 1} \frac{c_n}{n^s} &= \frac{1}{\Gamma (s)} \il^\infty_0 x^{s - 1} \left\{
 \sli^\infty_{n = 1} c_n (e^{- x})^n \right\} d x \, , \nonumber\\
  \sli^\infty_{n = 1}  c_n (e^{- x})^n &= \frac{1}{2 \pi i} \il^{a +  i \infty}_{a - i \infty} \frac{\Gamma 
  (s)}{x^s} \left\{ \sli^\infty_{n = 1} \frac{c_n}{n^s} \right\} d s \, .
\end{align}
One can see that the Mellin transform  changes the power series 
$\Sigma c_n (e^{-x})^n$  into a Dirichlet series $\Sigma c_n/n^s$
and the inverse of the Mellin transform changes the Dirichlet series into a power series.
\bi

\no
In particular, if we set $c_n = 1$ for all $n$, then with $\Sigma (e^{-x})^n  =  1/( e^{x}-1)$
we obtain an integral representation of the Riemann zeta function:        
\be
\zeta (s) = \sli^\infty_{n = 1} \frac{1}{n^s} = \frac{1}{\Gamma (s)} \il^\infty_0
\frac{x^{s - 1}}{e^x - 1} d x \, , \quad \quad Re (s) > 1 \, 
\ee
 the inverse of which is given by         
 \be
 \frac{1}{e^x - 1} = \frac{1}{2 \pi i} \il^{a + i \infty}_{a - i \infty} \frac{\Gamma (s) \zeta (s)}{x^x} d s \quad
 (a > 1) \, .
             \ee                                 
One of the most important formulae in Riemann's paper is given by
\be
\frac{\log \zeta (s)}{s} = \il^\infty_0 \Pi (x) x^{- s - 1} d x \, .
\ee
Here one recognizes for the first time the close connection between the zeta 
function and the function $\Pi (x)$. 
To understand the above formula better, let us take the logarithm of both sides of
\be
\zeta (s) = \pli_p \frac{1}{1 - p^{- s}}
\ee
and using $\log (1-x) = -x - 1/2 \, x^2  - 1/3 \,  x^3 \cdots$
we obtain
\be
\log \zeta (s) = - \sli_p \log (1 - p^{- s}) = \sum p^{- s} + \frac{1}{2} \sum p^{- 2s} + 
\frac{1}{3} \sum p^{- 3s} + \cdots \, .
\ee
Here we make use of the identities $(Re(s)>1)$
\be
p^{- s} = s \il^\infty_p x^{- s - 1} d s \, , \quad \quad p^{- 2s} = s \il^\infty_{p^2} x^{- s - 1} d x 
\, , \quad \quad \cdots \, , \quad \quad p^{-   ns} = s \il^\infty_{p^n} 
x^{- s - 1} d x \, , \cdots
\ee
to write
\begin{align}
 \log \zeta (s) &= \sli_p \sli_n \frac{1}{n} p^{- n s} = \sli_p \sli_n \frac{1}{n} \cdot s 
 \il^\infty_{p^n} x^{- s - 1} d x \nonumber\\
 &= s \il^\infty_0 \Pi (x) x^{- s - 1} d x \, .
\end{align}
To explain the last line, let us write
\begin{align}
 s \il^\infty_0 \Pi (x) x^{- s - 1} d x &= s \left\{ \left[ \cancel{\Pi (x) (- 1) \frac{1}{2} x^{- s}} \right]^\infty_0 - \il^\infty_0
 d x d \Pi \frac{x^{- s}}{- s} \right\} \nonumber\\
 &= \il^\infty_0 x^{- s} d \Pi (x) \quad \mbox{(Stieltjes integral)} \, ,
 \end{align}
where the measure $d \Pi$ has been written as the density times $dx$; more precisely:
\be
d \Pi = \lk \frac{d \Pi}{d x} \right) d x \, ,
\ee
where $d \Pi / dx$ is the density of primes plus $1/2$-density of prime squares,
plus $1/3$-density of prime cubes, etc.
\bi

\no
Let us not forget that the calculus version of the ``golden formula''
\be
\frac{\log \zeta (s)}{s} = \il^\infty_0 \Pi (x) x^{- s - 1} d x
\ee
has its origin in the Euler-Riemann prime product formula for the zeta function and the
intelligent invention of the step function  $\Pi(x)$. 
This name is justified because when $x$ is the exact square of a prime, e.g., 
$x = 9 = 3^2, \Pi (x)$
jumps up one-half, since $\pi (\sqrt{x}) = \pi (3)$ jumps up $1$, 
and so on. Note that the actual point where the jump occurs, the value of the function 
is halfway up the jump.
\bi

\no
So we have derived the marvelous formula given above, which will lead us directly to the
central result of Riemann's paper. But what is the inverted expression, i.e., how can we 
express $\Pi (x)$ in terms of  $\zeta (x)$? 
This will be discussed in the next chapter.
\bi

\no
\section{Riemann as an Expert in Fourier Transforms}
\bi

\no
Earlier we introduced the pair of equations
\begin{align}
 \frac{\log \zeta (s)}{s} &= \il^\infty_0 \Pi (x) x^{- s - 1} d x \quad \quad (Re (s)
 > 1) \, , \nonumber\\
 \mbox{and} \quad \Pi (x) &= \frac{1}{2 \pi i} \il^{a + i \infty}_{a - i \infty} 
 \log \zeta (s) x^s \frac{ds}{s} \quad \quad (a > 1) \, ,
\end{align}
when we discussed the Mellin transform. Let us see how Riemann reached the same result 
much earlier by employing the Fourier inversion formula:
\be
\varphi (x) = \frac{1}{2 \pi} \il^{+ \infty}_{- \infty}
\left[ \il^{+ \infty}_{- \infty} \varphi  (\lambda) e^{i (x - \lambda) \mu} d \lambda \right] d \mu \, .
\ee
When we write
\be
\varphi (x) = \il^{+ \infty}_{- \infty} \phi (\mu)  e^{i \mu x}
d \mu \, ,
\ee
we can consider $\phi (\mu)$  as coefficients of an expansion defined by 
\be
\phi (\mu) = \frac{1}{2 \pi} \il^{+ \infty}_{- \infty} \varphi (\lambda)
e^{- i \lambda \mu} d \lambda \, .
\ee
Now let $s = a+i \mu, a = const.>1$ and $\mu$ be a real variable.

Then with $\lambda = \log x$ and $\varphi (x) = 2 \Pi (e^x)e^{-ax}$, we obtain
\begin{align}
 {x = e^\lambda \atop \frac{d x}{x} = d \lambda} : \frac{\log \zeta (a + i \mu)}{a + i \mu} &= 
 \il^{+ \infty}_{- \infty}
 \Pi (e^\lambda) e^{- (a + i \mu) \lambda} d \lambda \, \nonumber\\
 =: \phi (\mu) &= \frac{1}{2 \pi} \il^{+ \infty}_{- \infty} \varphi (\lambda) e^{- i \mu \lambda}
 d \lambda \, .
 \end{align}
Hence we can continue to write
\be
(\varphi (x)) = 2 \pi \Pi (e^x) e^{- ax} = \il^{+ \infty}_{- \infty} 
\frac{\log \zeta (a + i \mu)}{a + i \mu} e^{i \mu x} d \mu
\ee
and using $e^x = y$, then $y \to x, s = a + i \mu, ds = id \mu, 
d\mu = 1/i \cdot d s$ 
we finally obtain
\be
\Pi (x) = \frac{1}{2 \pi i} \il^{a  + i \infty}_{a - i \infty} \log \zeta (s)
x^s \frac{ds}{s} \quad \quad (a > 1) \, ,
\ee
which is the desired result.

From here on we can directly arrive at Riemann's main result of his 1859 paper. 
However, for the time being we have to accept two of Riemann's novel quantities
(details will be reported later): The entire function $\xi (s)$ ($\zeta (s)$
is not an entire function) and the product formula for the $\xi$ function:
\begin{align}
 \xi (s) &= \frac{1}{2} s (s - 1) \pi^{- \frac{s}{2}} \Gamma 
 \lk \frac{s}{2} \right) \zeta (s) \, , 
 \quad \quad \Gamma \lk \frac{s}{2} \right) = \frac{2}{s} \Gamma \lk 1 + \frac{s}{2} \right)
 \nonumber\\
 &= (s - 1) \pi^{- \frac{s}{2}} \Gamma \lk 1 + \frac{s}{2} \right) \zeta (s)
\end{align}
and
\be
\xi (s) = \frac{1}{2} \pli_\rho \lk 1 - \frac{s}{\rho} \right) \, ,
\ee
with $\rho$ the zeros of the zeta function (equal to the zeros of $\xi$). 

So, taking the logarithm of both sides, we obtain
\begin{align}
 - \log 2 + \sli_p \log \lk 1 - \frac{s}{\rho} \right) &= \log (s - 1) - \frac{s}{2} \log
 \pi + \log \Gamma  \lk 1 + \frac{s}{2} \right) + \log \zeta (s) \nonumber\\
 \mbox{or} \quad \log \zeta (s) &= \sli_\rho \log \lk 1 - \frac{s}{\rho} \right) - \log 2 - \log \Gamma 
 \lk 1 + \frac{s}{2} \right) + \frac{s}{2} \log \pi - \log (s -  1) \, .
\end{align}
The first term on the right-hand side gives us the searched-for connection of the 
non-trivial zeta zeros with $\Pi (x)$. This becomes evident when we write
\be
\Pi (x) = \frac{1}{2 \pi i} \il^{a + i \infty}_{a - i \infty} \frac{\log \zeta (s)}{s} x^s
d s
\ee
with $\log \zeta (s)$ taken from above. Here, then, is Riemann's result:
\be
\boxed{\Pi (x) = Li (x) - \sli_\rho Li (x^\rho) + \log \lk \frac{1}{2} \right) + \il^\infty_x 
\frac{dt}{t  (t^2 - 1) \log t} \, , \quad \quad x > 1 \, .}
\ee
The sum over $\rho$ is to be understood as 
\be
\sli_{\,\mathrm{Im} \rho > 0}   (Li (x^\rho) + Li (x^{1 - \rho})) 
\ee
and $Li(x)$ denotes the logarithmic integral (see below).

This calculated expression for $\Pi (x)$ is then used in the formula
\be
\pi (x) = \sli^\infty_{n = 1} \frac{\mu (n)}{n} \Pi (x^{1/n}) = \Pi (x) - \frac{1}{2} \Pi
(x^{1/2}) - \frac{1}{3} \Pi
(x^{1/3}) - \frac{1}{5} \Pi
(x^{1/5}) + \frac{1}{6} \Pi
(x^{1/6}) + \cdots \, .
\ee
This is Riemann's great achievement, the explicit, exact calculation of the prime number counting function  $\pi(x)$.

Let us rewrite Riemann's result more explicitly:
\be
\Pi (x) = Li (x) - \sli_{\,\mathrm{Im} \rho > 0} (Li (x^\rho) + Li (x^{1 - \rho})) - \log 2 + \il^\infty_x \frac{dt}{t (t^2 - 1) \log t} \, , 
\quad \quad x > 1 
\ee
with
\be
Li (x) = \lim\limits_{\epsilon \to 0} \left[ \il^{1 - \epsilon}_{0} \frac{dt}{\log t} + \il^x_{1 + \epsilon} \frac{dt}{\log t} \right] \, .
\ee
If we differentiate $\Pi(x)$ we obtain
\be
d \Pi = \left[ \frac{1}{\log x} - \sli_{Re \alpha > 0} \frac{2 \cos (\alpha \log x)}{x^{1/2} \log x} - \frac{1}{x (x^2 - 1) \log x} \right]d x \, 
\quad \quad x > 1 \, .
\ee
$\alpha$ ranges  over all values such that $\rho = 1/2 + i \alpha$; in other words, $\alpha = -i(\rho - 1/2 )$ where $\rho$
ranges over all roots, so that 
\be
x^{\rho - 1} + x^{- \rho} = x^{- \frac{1}{2}} \left[ x^{i \alpha} + x^{- i \alpha} \right] = 2 x^{- \frac{1}{2}} \cos 
(\alpha \log x) \, .
\ee
The Riemann hypothesis says that the $\alpha$'s are all real.

Again, by the definition of $\Pi$, the measure $d \Pi$ is $dx$ times the density of primes plus $1/2$
the density of prime squares, plus $1/3$ the density of prime cubes plus, etc. Thus $1/(\log x)$
alone should not be considered an approximation only to the density of primes as Gau\ss{} suggested, but rather to $d\Pi/dx$, 
i.e., to the density of primes plus $1/2$ the density of prime squares, plus, etc.  

A fairly good approximation neglects the last term in $d\Pi$. It is the number of $\alpha$'s which is significant in $d\Pi$
which Riemann meant to study empirically to see the influence of the ``periodic terms''
on the distribution of primes. With the above equations we have reached the end of Riemann's famous paper of 1859.

We have, however, left out a number of revolutionary results to which we want to turn to now.

\section{On the Way to Riemann's Entire Function $\xi(s)$}

Let us begin with the integral representation of Euler's $\Gamma$ function:
\begin{align}
 \Gamma (s) &= \il^\infty_0 x^{s - 1} e^{- x} d x \, , \nonumber\\
 s \to \frac{s}{2}: \quad \Gamma \lk \frac{s}{2} \right) &= \il^\infty_0 x^{\frac{s}{2} - 1} e^{- x} d x \, , \nonumber\\
 x = \pi t n^2 : \quad \Gamma \lk \frac{s}{2} \right) &= \il^\infty_0  (\pi t n^2)^{\frac{s}{2} - 1} e^{- \pi t n^2}
 \pi n^2 dt \, , \nonumber\\
 \Gamma \lk \frac{s}{2} \right)  \pi^{- \frac{s}{2}} \frac{1}{n^s} &= \il^\infty_0  e^{- \pi t n^2} t^{\frac{s}{2}} \frac{dt}{t} \, , 
 \nonumber\\
 \mbox{Take} \quad \sli^\infty_{n = 1} : \quad \quad  \Gamma \lk \frac{s}{2} \right) \pi^{- \frac{s}{2}} \zeta (s) &= \il^\infty_0 \psi (t) t^{\frac{s}{2}} \frac{dt}{t} \, , 
 \quad \quad Re (s) > 1 \, , \nonumber\\
 \psi (t) &= \sli^\infty_{n = 1} e^{- \pi t n^2} \, .
\end{align}
The last equation defines one of Jacobi's $\vartheta$ functions:
\be
\Theta (x) := \vartheta_3 (0, i x) = \sli^{+ \infty}_{n = - \infty} e^{- \pi x n^2} \, , \quad \quad
\psi (x) = \sli^\infty_{n = 1} e^{- \pi x n^2} \, , \quad \quad 
\Theta (x) = 2 \psi (x) + 1 \, .
\ee
Also let me quote without proof the Jacobi identity:
\be
\Theta (x) = \frac{1}{\sqrt{x}} \Theta \lk \frac{1}{x} \right) \, , \quad \quad x > 0 \, .
\ee
One can then easily verify that
\be
\frac{1 + 2 \psi (x)}{1 + 2 \psi \lk \frac{1}{x} \right)} = \frac{1}{\sqrt{x}} \, , 
\ee
so that
\be
 \psi \lk \frac{1}{x} \right) = \frac{1}{2} \Theta  \lk \frac{1}{x} \right) - \frac{1}{2} = 
 \frac{1}{\sqrt{2}} \sqrt{x} \Theta (x) - \frac{1}{2} = \sqrt{x} \psi (x) + \frac{\sqrt{x}}{2} - \frac{1}{2} \, . \label{Gl: 72}
\ee
Now we are going to calculate the following integral, which will give us one of Riemann's wonderful results. 

Using $\Psi(x) = x^{-1/2} \Psi(1/x) - 1/2 + 1/2 x^{-1/2}$
and splitting the integral apart at 1, we obtain
\be
\il^\infty_0 \Psi (x) x^{s/2} \frac{dx}{x} = \il^\infty_1 \Psi (x) 
x^{s/2} \frac{dx}{x} + \il^1_0 \Psi \lk \frac{1}{x} \right) 
x^{\frac{s - 1}{2}} \frac{dx}{x} + \frac{1}{2} \il^1_0 
\lk x^{\frac{s-1}{2}}  - x^{\frac{s}{2}}
\right) \frac{dx}{x} \, .
\ee

In the last two integrals we substitute $x \to 1/x$ and so we get
\begin{align}
 \il^\infty_0 \Psi (x) x^{\frac{s}{2}} \frac{dx}{x} &= 
 \il^\infty_1 \Psi (x) \left[ x^{\frac{s}{2}} + x^{\frac{1}{2} (1 - s)} \right] 
 \frac{dx}{x} + \frac{1}{2} \il^\infty_1 \left[ x^{\frac{1}{2} (1 - s)} - x^{
 - \frac{s}{2}} \right] \frac{dx}{x} \nonumber\\
 \il^\infty_1 d x \left[ x^{- \frac{s}{2} - \frac{1}{2}} \right] &= - \frac{2}{s - 1}
 \, , \nonumber\\
 \il^\infty_1 d x \left[ x^{- \frac{s}{2} - 1} \right] &= \frac{2}{s} \, ,
 \nonumber\\
 &= \il^\infty_1 \Psi (x) \lk x^{\frac{s}{2} - 1} + x^{- \frac{s}{2} - \frac{1}{2}} 
 \right) d x = \frac{1}{s}  + \frac{1}{s - 1} \, .
\end{align}
Here, then, is the important formula contained in Riemann's paper:
\begin{align}
\Gamma \lk \frac{s}{2} \right) \pi^{- \frac{s}{2}} \zeta (s) = 
\il^\infty_1 \Psi (x) \lk x^{\frac{s}{2} - 1} + x^{- \frac{s}{2} - \frac{1}{2}}
\right) d x - & \frac{1}{s (1 - s)} \, .\nonumber\\
& pole \, \Gamma : s = 0 \atop 
pole \,  \zeta : s = 1
\end{align}
Notice that there is no change of the right-hand side under $s\to1-s$!
$\pi^{- s/2} \Gamma (s/2) \zeta (s)$ has simple poles at $s = 0$ and $s = 1$. To remove these poles, we multiply by
$1/2 s (s - 1)$. This is the reason why Riemann defines
\be
\xi (s) = \frac{1}{2} s (s - 1) \pi^{- \frac{s}{2}} \Gamma \lk \frac{s}{2} \right)
\zeta (s) \, ,
\ee
which is an entire function ($\zeta (s)$ is a meromorphic function.)
Obviously we have $\xi(s) = \xi(1-s)$
and the functional equation
\be
\Gamma \lk \frac{s}{2} \right) \pi^{- \frac{s}{2}} \zeta (s) = \Gamma
\lk \frac{1 - s}{2} \right) \pi^{- \frac{1}{2} (1 - s)} \zeta (1 - s) \, .
\ee
We obtain the right-hand side by the left-hand side by replacing 
$s$ by $(1-s)$.

Now we can continue to write for $\xi(s)$
\begin{align}
 \xi (s) &= \frac{1}{2} - \frac{s (1 - s)}{2} \il^\infty_1 \Psi (x) \lk x^{\frac{s}{2}} + x^{\frac{1}{2} (1 - s)} \right)
 \frac{dx}{x} \nonumber\\
 &= \frac{1}{2} - \frac{s (1 - s)}{2} \il^\infty_1  \frac{d}{dx} \left\{ \Psi (x) \left[ \frac{x^{\frac{s}{2}}}{\frac{s}{2}} + 
 \frac{x^{\frac{1}{2} (1 - s)}}{\frac{1}{2} (1 -  s)} \right] \right\} dx \nonumber\\
  & + \frac{s (1 - s)}{2} \il^\infty_1 \Psi' (x) \left[ \frac{x^{\frac{s}{2}}}{\frac{s}{2}} + \frac{x^{\frac{1}{2} (1 - s)}}{\frac{1}{2}
  (1 - s)} \right] ds \nonumber\\
  &= \frac{1}{2} + \frac{s (1 - s)}{2} \Psi (1) \left[ \frac{2}{3} + \frac{2}{1 - s} \right] \nonumber\\
  & + \il^\infty_1 \Psi' (x) \left[ (1 - s) x^{\frac{s}{2}} + s x^{\frac{1}{2} (1  -s )} \right] dx \nonumber\\
  &= \frac{1}{2} + \Psi (1) + \il^\infty_1  x^{\frac{s}{2}} \Psi' (x) \left[ (1 - s)   x^{\frac{1}{2} (s - 1) - 1} 
  + s x^{- \frac{s}{2} - 1} \right] dx \nonumber\\
  &= \frac{1}{2} + \Psi (1) + \il^\infty_1 \frac{d}{dx} \left[ x^{\frac{3}{2}} \Psi' (x) \lk -  2 x^{\frac{1}{2} (s - 1)} - 2 x^{- \frac{s}{2}} 
  \right) \right] dx \nonumber\\
  & - \il^\infty_1 \frac{d}{dx} \left[ x^{\frac{3}{2}} \Psi' (x) \right] \left[ - 2 x^{\frac{1}{2} (s - 1)} - 2 x^{- \frac{s}{2}} \right] dx 
  \nonumber\\
  &= \frac{1}{2} + \Psi (1) - \Psi' (1) [- 2 - 2] + \il^\infty_1 \frac{d}{dx} \left[ 
  x^{\frac{3}{2}} \Psi' (x) \right] 
  \lk 2 x^{\frac{1}{2} (s - 1)}  + 2 x^{- \frac{s}{2}} \right) d x \, .
\end{align}
Differentiation of 
\be
2 \Psi (x) + 1 = x^{- \frac{1}{2}} \left[ 2 \Psi \lk \frac{1}{x} \right) + 1 \right]
\ee
easily gives
\be
\frac{1}{2} + \Psi (1) + 4 \Psi' (1) = 0
\ee
and using this puts the formula in the final form:
\be
\xi (s)  = 4 \il^\infty_1 \frac{d}{dx} \left[ x^{\frac{3}{2}} \Psi' (x) \right] x^{- \frac{1}{4}} \cosh \left[ \frac{1}{2}   \lk s - \frac{1}{2}
\right) \log x \right] d x \, ,
\ee
or, as Riemann writes it $(s = 1/2+it; 1/2$ is Riemann's conjecture!):
\be
\Xi (t) = \xi \lk \frac{1}{2} + it \right) = 4 \il^\infty_1 \frac{d}{dx} \left[ x^{\frac{3}{2}} \psi' (x) \right] x^{- \frac{1}{4}}
\cos \lk \frac{t}{2} \log x \right) dx \, .
\ee
With
\be
\frac{d}{dx} \left[ x^{3/2} \psi' (x) \right] = \sli^{\infty}_{n = 1} \lk n^4 \pi^2 x - \frac{3}{2}
n^2 \pi \right) x^{1/2} \exp (-n^2 \pi x)
\ee
and
\be
v = \frac{1}{2} \log x
\ee
and then $v = 2 u$, we can also write $\Xi \lk \frac{t}{2} \right)$ as a Fourier transform
\be
\Xi \lk \frac{t}{2} \right) = 8 \il^\infty_{0} d u \Phi (u) \cos (ut)
\ee
with
\be
\Phi (u) = \sli^\infty_{n = 1} \pi n^2 \lk 2 n^2 \pi \exp (4 u) - 3 \right) \exp (5u - n^2  \pi \exp (4u)) \, .
\ee
If $cosh [1/2(s-1/2)\log x]$ is expanded in the usual power series
\be
\cos h y = \frac{1}{2} \lk e^{y} + e^{- y} \right) = \sli \frac{y^{2n}}{(2 n)!} \, ,
\ee
we can write
\be
\xi (s) = \sli^\infty_{n = 0} a_{2n} \lk s - \frac{1}{2} \right)^{2n} \, ,
\ee
where
 \be
 a_{2n} = 4 \il^\infty_1 \frac{d}{dx} \left[ x^{3/2} \Psi' (x) \right]  x^{- \frac{1}{4}} \frac{\lk \frac{1}{2}
 \log x \right)^{2n}}{(2 n)!} d x \, .
 \ee
Let us return to 
\be
\xi (s) = \frac{1}{2} s (s - 1) \pi^{- \frac{s}{2}} \Gamma \lk \frac{s}{2} \right) \zeta (s) \, ,
\ee
with
\be
\pi^{- \frac{s}{2}} \Gamma \lk \frac{s}{2} \right) \zeta (s) = \frac{1}{s (s - 1)} + \il^\infty_1 \Psi (x) 
\lk x^{\frac{s}{2} - 1} + x^{- \frac{s}{2} - \frac{1}{2}} \right) dx \, ,
\ee
and write the right-hand side in terms of $s =1/2 + it$, 
which makes use of Riemann's conjecture $Re(s) = 1/2$. 
Since the details of the substitution are trivial, we merely give the final result:
\begin{align}
 \xi \lk \frac{1}{2} + i t \right) &= \frac{1}{2} \lk \frac{1}{2} + i t \right) \lk it - \frac{1}{2}\right) \pi^{- \frac{1}{4} - 
 i \frac{t}{2}} \Gamma \lk \frac{1}{4}  + i \frac{t}{2} \right) \zeta \lk \frac{1}{2} + i t \right) \nonumber\\
 &= \frac{- \lk t^2  + \frac{1}{4} \right)}{\left[ 2 (\sqrt{\pi})^{\frac{1}{2} + it} \right]} \Gamma \lk 
 \frac{1}{4} + \frac{it}{2} \right) \zeta \lk \frac{1}{2} + it \right) \, .
 \end{align}
 In particular,
 \be
 \xi \lk \frac{1}{2} \right) = \frac{- 1}{(8 \pi^{1/4})} \Gamma \lk \frac{1}{4} \right) \zeta \lk \frac{1}{2} \right)
 \ee
 with
 \be
 \zeta \lk \frac{1}{2} \right) = - 1.4603545088 \, , \quad \Gamma \lk \frac{1}{4} \right) = \sqrt{2 \varpi 2 \pi} = 
 3.6256099082 \, ,
 \ee
 where Gauss' lemniscate constant is given by 
 \be
 \varpi = 2.62205755429 \, .
 \ee
 Altogether:
 \be
 \xi \lk \frac{1}{2} \right) = 0.4971207781 = a_0 \, , 
 \ee
 which is the minimum for the real valued $\xi (s)$ at $s = 1/2$.
 By the way $\xi (0) = \xi (1) = - \zeta (0) = 1/2$. The above result can also be written as
 \begin{align}
 \Xi (t) := \xi \lk \frac{1}{2} + i t \right) &= \frac{1}{2} - \lk t^2 + \frac{1}{4} \right) \il^\infty_1 
 \Psi (x) x^{- \frac{3}{4}} \cos \lk \frac{t}{2} \log x \right) dx \, .
\end{align}
The right-hand side of this equation tells us that because $t \in R_e, x \in R_e$ and $\log x \in R_e$, we have
\be
\,\mathrm{Im} \xi \lk \frac{1}{2} + i t \right) = 0 \, , \quad \quad i.e., \quad \xi \lk \frac{1}{2} + i t \right) \equiv \Xi (t) \in R_e \, . \label{Gl: 98}
\ee
Since  $\Xi (t) = \xi(1/2+it)$ for $t \to \infty$ changes its sign infinitely often, $\xi(s)$ (and $\zeta (s)$) 
must have infinitely many zeros on $Re(s) = 1/2$.

There is another useful form $\xi(s)$ that starts with its original definition:
\begin{align}
 \xi (s) &= \frac{s (s - 1)}{2} \Gamma \lk \frac{s}{2} \right) \pi^{- \frac{s}{2}} \zeta (s) \nonumber\\
 &= e^{\ln \Gamma  \lk \frac{s}{2} \right) } \pi^{- \frac{s}{2}} \frac{s (s - 1)}{2}  \zeta (s) \, . \label{Gl: 99}
\end{align}
Then, setting $s = 1/2 + it$, we have
\begin{align}
 \xi  \lk \frac{1}{2} + it  \right) &= e^{\ln \Gamma \frac{\lk \frac{1}{2} + i t \right)}{2}} \pi^{- \frac{\frac{1}{2}  + it}{2}} \frac{1}{2} 
 \lk \frac{1}{2} + it \right) \lk \frac{1}{2} + it - 1 \right) \zeta \lk \frac{1}{2} + it \right) \nonumber\\
 &= \left[ e^{R_e \ln \frac{\lk \frac{1}{2} + it \right)}{2}} \pi^{- \frac{1}{4}} \cdot \frac{- t^2 - \frac{1}{4}}{2}
 \right] \left[ e^{i \,\mathrm{Im} \ln \Gamma \frac{\lk \frac{1}{2} + it \right)}{2} } \pi^{- \frac{it}{2}} \zeta 
 \lk \frac{1}{2} + it \right) \right] \nonumber\\
 &= \left[ - e^{R_e \ln \Gamma \lk \frac{\frac{1}{2} + it}{2} \right)} \pi^{- \frac{1}{4}} \frac{t^2 + \frac{1}{4}}{2} \right]
 \left[ e^{i \,\mathrm{Im} \ln \Gamma \lk \frac{\frac{1}{2} + it}{2} \right)} \pi^{- \frac{it}{2}} \zeta \lk \frac{1}{2} + it \right) \right] \, .
\end{align}
Notice that the first factor in the square brackets is negative. For the second factor we have 
\be
Z (t) = e^{i \vartheta (t)} \zeta \lk \frac{1}{2} + it \right) \, , \quad \quad
\vartheta (t) = \,\mathrm{Im} \ln \Gamma \lk \frac{\frac{1}{2} + i t}{2} \right) - \frac{t}{2} \ln \pi \, .
\ee
Thus, $Z(t)$ has always the opposite sign compared to the $\xi$ function.

Now we have to compute $\vartheta (t)$ and $\zeta(1/2 + it)$. For numerical analysis it is sufficient to use
\be
\vartheta (t) \sim \frac{t}{2} \log \frac{t}{2 \pi} - \frac{t}{2} - \frac{\pi}{8} + \frac{1}{48  t} \, ,
\ee
which one can then apply to compute the roots of $\xi(s)$ on the critical line. 
\bi

\no
\section{The Product Representation of  $\xi(s)$ and $\zeta (s)$ by Riemann (1859) and Hadamard (1893)}
\bi

\no
Riemann's goal (before Weierstrass!) was to prove that $\xi(s)$ can be expanded as an infinite product
\be
\xi (s) = \xi (0) \pli_\rho \lk 1 - \frac{s}{\rho} \right) \, ,
\ee
where $\rho$ ranges over all the roots of $\xi(\rho) = 0$. 
He did not really prove this formula, but he was right, as shown much later by Hadamard. But one has to admit that Riemann must 
have had a strong inkling of the product formula Weierstrass was soon to introduce as an essential representation into the function 
theory, here the entire functions, i.e., functions that can be determined  by their zeros. 

As a brief reminder, here is Weierstrass' product representation of the $\Gamma$ function:
\be
\Gamma (x) = e^{- \gamma x} \frac{1}{x} \pli^\infty_{k = 1} \frac{e^{\frac{x}{k}}}{\lk 1 + \frac{x}{k} \right)} \, ,
\ee
where $\gamma$ is the Euler-Mascheroni constant,
\be
\gamma = \lim\limits_{n \to \infty} \left[ \sli^n_{k = 1} \frac{1}{k} - \log n \right] \simeq 0.5772157 \, .
\ee
From this product formula follows, with the aid of
\be
\Gamma (x) \Gamma (1 - x) = \Gamma (x) (- x) \Gamma (- x) = \frac{\pi}{\sin (\pi x)} \, , 
\ee
the product representation of $\sin (\pi x)$:
\begin{align}
 \sin (\pi x) &= - \frac{\pi}{x} \frac{1}{\Gamma (x) \Gamma (- x)} = - \frac{\pi}{x} 
 \lk e^{\gamma x} x \pli^\infty_{k = 1} \frac{\lk 1 + \frac{x}{k} \right)}{e^{\frac{x}{n}}} \right)
 \lk e^{- \gamma x} (- x) \pli^\infty_{k > 1} \frac{\lk 1 - \frac{x}{k} \right)}{e^{- \frac{x}{k}}} \right)\nonumber\\
 &= \pi x \pli^\infty_{k = 1} \lk 1 - \frac{x^2}{k^2} \right) \, ,
\end{align}
a polynomial of infinite degree. Similarly, Euler thought of  $\sin (\pi x)$
as a polynomial of infinite degree when he conjectured, and finally proved, the formula for $\sin (\pi x)$.

So, why not think of $\xi(s)$
as a polynomial of infinite degree and write down a product formula determined by its infinite zeros $\rho$? 
This is what Hadamard had done in 1893 in a paper in which he studied entire functions and their representations as infinite products 
-- like Weierstrass. He was able to prove that Riemann's product formula was correct:
\be
\xi (s) = \xi (0) \pli_\rho \lk 1 - \frac{\xi}{\rho} \right) \, .
\ee
$\xi(s)$ is an entire function. The infinite product is understood to be taken in an order which pairs each root $\rho$
with the corresponding root $1-\rho$.  Hadamard's proof of the product formula for $\xi$
was called by von Mangoldt ``the first real progress in the field in 34 years,'' that is, the first since Riemann.

Hadamard showed that it is possible to construct the $\zeta$ function as an infinite product, given its zeros:
\be
\zeta (s) = f (s) \pli_{\zeta (\rho) = 0} \lk 1 - \frac{s}{\rho} \right) e^{\frac{s}{\rho}} \, , 
\quad \quad f (s) = \frac{1}{2 (s - 1)} \lk \frac{ 2 \pi}{e}
\right)^s \, .
\ee
Hence, including the trivial as well as the non-trivial zeros he obtains
\be
\zeta (s) = \frac{1}{2 (s - 1)} \lk \frac{2 \pi}{e} \right)^s \pli^\infty_{n = 1} \lk 1 + 
\frac{s}{2n} \right) e^{- \frac{s}{2n}} \cdot \pli_\rho \lk 1 - \frac{s}{\rho} \right) e^{\frac{s}{\rho}} \, .
\ee
For the first product we use the product representation given by Weierstrass:
\be
\frac{1}{\Gamma (1 + s)} = e^{\gamma s} \pli^\infty_{n = 1} \lk 1 + \frac{s}{n} \right) e^{- \frac{s}{n}} \, ,
\ee
and so obtain the Hadamard product formula, which is convergent in $C  \setminus \{1\}$:
\be
\zeta (s) = \frac{e^{\lk \log 2 \pi - 1 - \frac{\gamma}{2} \right) s}}{2 (s - 1) \Gamma \lk 
1 + \frac{s}{2} \right) } \pli_\rho \lk 1 - \frac{s}{\rho} \right) e^{\frac{s}{\rho}} \, .
\ee
A slightly simplified form of the Hadamard product is
\be
 \zeta (s) = \frac{\pi^{s/2}}{2 (s - 1) \Gamma \lk 1 + \frac{s}{2} \right)} \pli_\rho \lk 1 - \frac{s}{\rho} \right) \, .
\ee
Here we took pairs of roots $\rho$ and $-\rho$ together so the exponents $e^{-s / \rho}$ cancel.

The last expression shows the the $\zeta$ function can be completely constructed by its roots (Riemann's specialty) and the 
singularity at $s =1$. However, to obtain absolute convergence, we have to introduce $\rho$ and $-\rho$ pairwise 
in the product.

Now, we remember Riemann's entire function $\xi(s)$ and how it is related to the (non-entire) $\zeta$ function:
\be
\xi (s) = \frac{s (s - 1)}{2} \Gamma \lk \frac{s}{2} \right) \pi^{- \frac{s}{2}} \zeta (s) \, .
\ee
Then
\be
\frac{s (s - 1)}{2} \pi^{- \frac{s}{2}} \Gamma \lk \frac{s}{2} \right) \cdot
\frac{\pi^{s/2}}{2 (s - 1) \Gamma \lk 1 + \frac{s}{2} \right)} \pli_\rho \lk 1 - \frac{s}{\rho} \right) \, , \quad \quad 
\Gamma \lk 1 + \frac{s}{2} \right) = \frac{s}{2} \Gamma \lk \frac{s}{2} \right)
\ee
or
\be
\xi (s) = \frac{1}{2} \pli_\rho \lk 1 - \frac{s}{\rho} \right)
\ee
and using $\xi (0) = \frac{1}{2}$, we have indeed
\be
\xi (s) = \xi (0) \pli_\rho \lk 1 - \frac{s}{\rho} \right) \, ,
\ee
 which is Riemann's result of 1859! 
 
 Since the zeros of $\zeta (s)$ and $\xi(s)$ in the critical strip are identical, we can also write
 \begin{align}
  \zeta (s) &= \frac{\pi^{s/2}}{2 (s - 1) \Gamma \lk 1 + \frac{s}{2} \right) } \pli_\rho
  \lk 1 - \frac{s}{\rho} \right)    \lk 1 - \frac{s}{1 - \rho} \right) \nonumber\\
  &= \frac{\pi^{s/2}}{2 (s - 1) \Gamma \lk 1 + \frac{s}{2} \right) }  \lk 1 -
  \frac{s}{\frac{1}{2} + 14.134 i} \right) 
  \lk 1 -
  \frac{s}{\frac{1}{2} - 14.134 i} \right) \lk 1 -
  \frac{s}{\frac{1}{2} + 21.022 i} \right)  \lk \cdots \right) \, ,
 \end{align}
where we have used the first zeros on the $Re(s) =1/2$ axis.
\bi

\no
\section{Derivation of Von Mangoldt's Formula for  $\Psi (x)$}
\bi

\no
There is another, more modern version of an equivalent to Riemann's formula for $\Pi(x)$, i.e., 
\be
\Pi (x) = Li (x) - \sli_\rho Li (x^\rho) + \log \xi (0) + \il^\infty_x \frac{dt}{t (t^2 - 1) \log t} \quad (x > 1) \, .
\ee
This is von Mangoldt's formula for $\Psi (x)$, which contains essentially the same information as Riemann's $\Pi(x)$. 
On the way to the explicit formula for $\Psi (x)$, we need a special representation of the discontinuity function. 
So let us begin very simply by verifying
\begin{align}
 \frac{1}{s - \beta} &= \il^\infty_1 x^{- s} x^{\beta - 1} dx \, , \quad \quad Re (s - \beta) > 0 \, , \nonumber\\
 x = e^\lambda: &= \il^\infty_0 e^{- \lambda s} e^{\lambda (\beta - 1)} e^\lambda d \lambda = \il^\infty_0 
 e^{- \lambda s} e^{\lambda \beta} d \lambda \, , \nonumber\\
 s = a + i \mu &= \il^\infty_0 e^{- \lambda (a + i \mu)} e^{\lambda \beta} \alpha \lambda \, , \nonumber\\
 \frac{1}{a + i \mu - \beta} &= \il^\infty_0 e^{- i \lambda \mu} e^{\lambda  (\beta - a)} d \lambda \, , 
 \quad \quad a > Re \beta \, , \nonumber\\
 \il^{+ \infty}_{- \infty} \frac{1}{a +   i \mu - \beta} e^{i \mu x} d \mu &= \il^{+ \infty}_{- \infty} e^{i \mu x} d \mu \il^\infty_0 
 e^{- i \lambda \mu} e^{\lambda (\beta - a)} d \lambda 
 \nonumber\\
 &= \il^{+ \infty}_{- \infty} \left[ \il^\infty_0 e^{i (x - \lambda) \mu} d \mu \right] e^{\lambda (\beta - a)} d \lambda \nonumber\\
 &= \il^{+ \infty}_{- \infty} 2 \pi \delta (x - \lambda) e^{\lambda (\beta - a)} d \lambda \nonumber\\
 &= \left\{ \begin{array}{ccc}
             2 \pi e^{x (\beta- a)} & , & x > 0 \\
             0 & , & x < 0 
            \end{array} \right. \, .
\end{align}
So far we have
\be
\frac{1}{2 \pi} \il^{+ \infty}_{- \infty} \frac{1}{a + i \mu - \beta} e^{x (a + i \mu)} d \mu = 
\left\{ \begin{array}{ccc}
        e^{x \beta} & , & x > 0 \\
        0 & , & x < 0 
       \end{array}
       \right. \, . \ee
With $e^x = y$ and $s = a + i \mu$, we obtain the discontinuity factor (step function) 
\begin{align}
 \frac{1}{2 \pi i} \il^{a + i \infty}_{a - i \infty} \frac{1}{s - \beta} y^s d s &=  \left\{
 \begin{array}{ccc}
  y^\beta & , & y > 1 \\
  0 & , & y < 1
 \end{array}
 \right. 
 \stackrel{\beta = 0}{=}
 \left\{ 
 \begin{array}{ccc}
  1 & , & y > 1 \\
  \frac{1}{2} & , & y = 0 \\
  0 & , & y < 1
 \end{array}
 \right. 
 \quad \quad a > 0 
 \, .
\end{align}
Now we go back to the Euler-Riemann zeta function,
\be
\zeta (z) = \pli_{p \in P} \frac{1}{1 - p^{- z}} \, , \quad \quad Re (z) > 1 
\ee
and take the logarithm:
\begin{align}
 \log \zeta (z) &= - \sli_p \log (1 - p^{- z}) = - \sli_p \log \lk 1 - e^{- z \log p} \right) \, ,   
 \nonumber\\
 \frac{d}{d z} \log \zeta (z) &= - \sli_p \frac{1}{1 - p^{-  z}} \frac{d}{d z} \lk 1 - e^{ - z \log p} \right) 
 = - \sli_p \frac{1}{1 - p^{- z}} \log p \cdot p^{- z} \nonumber\\
 &= - \sli_p \frac{p^{- z}}{1  - p^{- z}} \log p = - \sli_p \sli^\infty_{\nu = 1} p^{- \nu z} \log p \nonumber\\
 &= \frac{\zeta' (z)}{\zeta (z)}  \, . \nonumber\\
 \cdot \frac{x^z}{z} \, : \, \frac{x^z}{z} \sli^\infty_{p  \atop \nu = 1} \frac{\log p}{p^{\nu z}} &=
 \sli^\infty_{p \atop \nu = 1} \lk \frac{x}{p^\nu} \right) \frac{\log p}{z} = - \frac{\zeta' (z)}{\zeta (z)} \cdot \frac{x^z}{z} \, , \nonumber\\
 \frac{1}{2 \pi i} \il^{a + i \infty}_{a - i \infty} \sli^\infty_{p, \nu = 1} \lk \frac{x}{p^\nu} \right)^z \frac{\log p}{z} &=
 \frac{1}{2 \pi i} \il^{a + i \infty}_{a - i \infty} - \frac{\zeta' (x)}{\zeta (z)} \frac{x^z}{z} d z \nonumber\\
 \mbox{or} \quad \sli^\infty_{p \atop \nu = 1} \log p \frac{1}{2 \pi i} \il^{a + i \infty}_{a - i \infty} 
 \lk \frac{x}{p^\nu} \right)^z \frac{1}{z} d z &=
 \frac{1}{2 \pi i} \il^{a + i \infty}_{a - i \infty} - \frac{\zeta' (z)}{\zeta (z)} \frac{x^z}{z} d z \nonumber\\
 y = \frac{x}{p^\nu} : \quad \sli^\infty_{p \atop \nu = 1} \log p \frac{1}{2 \pi i} \il^{a + i \infty}_{a - i \infty}
 \frac{y^z}{z} dz &= \frac{1}{2 \pi i} \il^{a + i \infty}_{a - i \infty} - \frac{\zeta' (z)}{\zeta (z)} \frac{x^z}{z} d z \, .
\end{align}
Here we use the $1$ of the discontinuity factor on the left-hand side and so obtain the Chebyshev function $\Psi (x)$:
\be
\Psi (x) = \sli_{p^\nu < x} \log p = \frac{1}{2 \pi i} \il^{a + i \infty}_{a - i \infty} - \frac{\zeta' (z)}{\zeta (z)} \frac{x^z}{z} d z \, .
\ee
So one has to sum the logarithm of all primes up to $x$. $p^\nu > x$ would mean $y<1$, 
but for this case the discontinuity formula gives zero.

The integral of the right-hand side can be evaluated with the aid of the theorem of residues. 
The contributions to the residues of $\zeta'(z)/\zeta(z) \cdot x^z/z$ come from
\begin{align}
 \begin{array}{ccc}
  \mbox{Singularity} & \mbox{Reason} & \mbox{Residue} \\
  0 & \frac{x^z}{z} & \frac{\zeta' (0)}{\zeta (0)} = \frac{- \frac{1}{2} \log 2 \pi}{- \frac{1}{2}} = \log (2 \pi) \\
  1 & \mbox{pole of} \, \zeta \, \quad \quad \frac{\zeta' (z)}{\zeta (z)} = - \frac{1}{z - 1} + \gamma + \cdots & 
  \lim\limits_{z \to 1} (z - 1) \lk \frac{- 1}{z - 1} + \mathcal{O} (1) \right) \frac{x^z}{z} = \frac{- x^1}{1} = - x \\
  - 2, - 4, - 6, \cdots & \mbox{trivial zeros of} \, \zeta (z) & \frac{1}{2} x^{- 2}, \frac{1}{4} x^{- 4}, \frac{1}{6}
  x^{- 6} , \cdots \atop \sli^\infty_{n = 1} \frac{x^{- 2n}}{2 n}  = \frac{1}{2} \log  \lk 1 - \frac{1}{x^2} \right) \\
  \rho & \mbox{nontrivial zeros of} \, \zeta (z) & \frac{x^\rho}{\rho}
 \end{array}
\end{align}
which leads to the exact explicit formula
\be
\Psi (x) =  x - \log (2 \pi) - \frac{1}{2} \log\lk 1 - \frac{1}{x^2} \right) - \sli_{\zeta (\rho) = 0} \frac{x^\rho}{\rho} \, .
\ee
This is known as Mangoldt's formula (1895) and is one of the most important formulae in analytic theory of numbers. 
$\Psi (x)$ is real and gives the jumps for prime powers $x$. 
Although the last term looks complex, it is not, since the zeros enter pairwise and hence it is also real.

$\Psi (x)$ is equivalent to Riemann's $\Pi(x)$ and one has to admit that the formula for $\Psi (x)$ 
was deduced much more easily than the formula for $\Pi(x)$, with which we began this chapter. 
No wonder that it is meanwhile considered preferable to that of $\Pi(x)$.
\bi

\no
\section{The Number of Roots in the Critical Strip}

The following theorem was originally formulated by Riemann -- but not proved. It was not until 1905 that von Mangoldt proved that
the number of zeros of $\zeta$ in the critical range $0<Re(s)<1, 0<t<T$ is given by
\be
N (T) = \frac{T}{2 \pi} \log \frac{T}{2 \pi} - \frac{T}{2 \pi} \, .
\ee
To prove this statement, let us assume $T \geq 3$ and $\zeta(s) \neq 0$ for $t = T$.
\bi

\no
Then consider the rectangular $R_T$ in the complex plane:
\begin{figure}[h]
\centerline{  
\includegraphics[width=6cm]{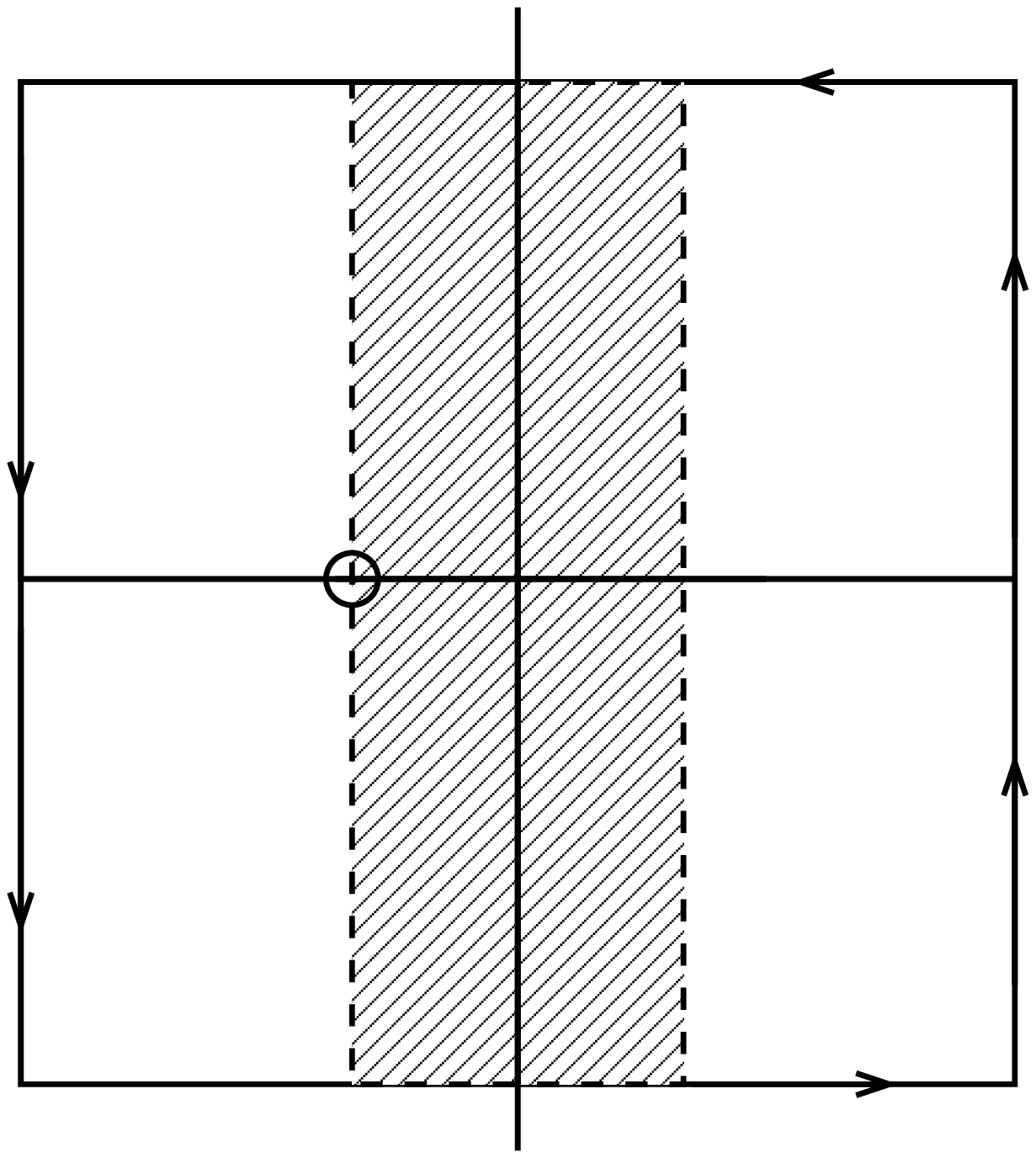}
\put(-5,0){$2 - iT$}
\put(-105,-5){crit. strip}
\put(-205,0){$-1 - iT$}
\put(5,60){$C_T$}
\put(-190,50){$C'_T$}
\put(0,80){$2$}
\put(-55,85){$+ 1$}
\put(-185,85){$-1$}
\put(-80,103){$\frac{1}{2}$}
\put(-105,103){$\frac{1}{2}$}
\put(5,175){$2 + iT$}
\put(-210,175){$-1 + iT$}
} 
\caption{Boundary of $R_T$}\label{plane}
\end{figure}

The zeros of the $\xi$ function are identical to the ones of the $\zeta$
function in the critical range. Symmetry with respect to the axis $Re(s) =  1/2$ yields (remember from the logarithmic residue)   
\be
2 N (T)  = \frac{1}{2 \pi i} \il_{\partial R_T} \frac{\xi' (s)}{\xi (s)} d s \, .
\ee
From the functional equation of $\xi$ we obtain
\begin{align}
 \xi (1 - s) &= \xi (s) \nonumber\\
 - \frac{\xi' (1 - s)}{\xi (1 - s)} &= \frac{\xi' (s)}{\xi (s)} \, .
\end{align}
$C'_T (C_T)$ is the left (right) boundary of $R_T$:
\begin{align}
 \il_{C'_T} \frac{\xi' (s)}{\xi (s)} d s &= \il_{C_T} \frac{\xi' (1 - s)}{\xi (1 - s)} d (1 - s) = \il_{C_T} \frac{\xi' (s)}{\xi (s)} ds \nonumber\\
 > \quad \quad N (T) &= \frac{1}{2 \pi i} \il_{C_T} \frac{\xi' (s)}{\xi (s)} d s  \, .
\end{align}
Now, using the following representation of the $\xi$ function,
\be
\xi (s) = \frac{s (s - 1)}{2} \pi^{- s/2} \Gamma \lk \frac{s}{2} \right) \zeta (s)
\ee
we take the logarithm
\begin{align}
 \log \xi (s) &= - \log 2 + \log s + \log (s - 1) - \frac{s}{2} \log \pi + \log \Gamma \lk 
 \frac{s}{2} \right) + \log \zeta 
 \nonumber\\
 > \frac{d}{ds} \log \xi (s) &= \frac{\xi' (s)}{\xi (s)}  = \frac{1}{s} + \frac{1}{s - 1} - \frac{1}{2} \log \pi + 
 \frac{1}{2} \frac{\Gamma' \lk \frac{s}{2} \right)}{\Gamma \lk \frac{s}{2} \right)} + \frac{\zeta' (s)}{\zeta (s)}
  \nonumber\\
  > 2 \pi i N (T) &= \underbrace{\il_{C_T} \lk \frac{1}{s} + \frac{1}{s - 1} \right) 
  ds}_{\textcircled{\raisebox{-0.3mm}{1}}} - 
  \underbrace{\il_{C_T} \frac{1}{2} \log \pi 
  ds}_{\textcircled{\raisebox{-0.3mm}{2}}} + 
  \underbrace{\frac{1}{2} \il_{C_T} \frac{\Gamma' \lk \frac{s}{2} \right)}{\Gamma \lk \frac{s}{2} \right)}\,
  ds}_{\textcircled{\raisebox{-0.3mm}{3}}} + \il_{C_T} 
  \frac{\zeta' (s)}{\zeta (s)} d s 
\end{align}
\begin{align}
  \textcircled{\raisebox{-0.3mm}{1}} &\quad& \il_{C_T} \lk \frac{1}{s} + \frac{1}{s - 1} \right) d s &= \frac{1}{2} 
  \il_{\partial R_T} \lk \frac{1}{s} + \frac{1}{s - 1} \right)
  d s \stackrel{\mbox{resid}}{=} \frac{1}{2} 2 \pi i (1 + 1) = 2 \pi i 
\nonumber\\
  \textcircled{\raisebox{-0.3mm}{2}} &\quad& \il_{C_T} \frac{1}{2} \log \pi d s &= \frac{1}{2} \log \pi \lk \lk \frac{1}{2} + i T \right) - 
  \lk \frac{1}{2} - i T \right) \right) = i T \log \pi 
\nonumber\\
 \textcircled{\raisebox{-0.3mm}{3}} &\quad& \il_{C_T} \frac{1}{2} \frac{\Gamma' \lk \frac{s}{2} \right)}{\Gamma \lk \frac{s}{2} \right)} d s &= \left. \log
 \Gamma \lk \frac{s}{2} \right) \right|^{\frac{1}{2} + i  T}_{\frac{1}{2} - i T} 
\nonumber\\
&& &= \log \Gamma \lk \frac{1}{4} + i \frac{T}{2} \right) - \log \Gamma \lk \frac{1}{4} - i \frac{T}{2} \right)  
\end{align}
\begin{align}
 \log \Gamma (\bar{s}) = \overline{\log \Gamma (s)}: &= 2 i \,\mathrm{Im} \log \Gamma \lk \frac{1}{4} + i \frac{T}{2} \right)  
\nonumber\\
 \overset{\text{Expand}}{\underset{T \geq 3}{=}}
 &= 2 i \,\mathrm{Im} \lk \log \sqrt{2 \pi} + \lk - \frac{1}{4} + i \frac{T}{2} \right) \log 
 \lk i \frac{T}{2} \right)  - i \frac{T}{2}  + \mathcal{O}  \lk \frac{1}{T} \right) \right) 
\nonumber\\
 &= 2 i \,\mathrm{Im} \lk \log \sqrt{2 \pi} + \lk - \frac{1}{4} + i \frac{T}{2} \right)  \lk \log \frac{T}{2}  + i \frac{\pi}{2} \right) - i \frac{T}{2}
 + \mathcal{O}\left( \frac{1}{T}\right) \right)  
 \nonumber\\
  &= 2 i \pi \lk \frac{T}{2 \pi} \log \frac{T}{2} - \frac{T}{2 \pi} \right) - \frac{1}{8} + \mathcal{O} \lk \frac{1}{T} \right) \, .
 \end{align} 
Our intermediate result is then
\be
2 \pi i N (T) = 2 \pi i - i T \log \pi + 2 \pi i \lk 
\frac{T}{2 \pi} \log \frac{T}{2} - \frac{T}{2 \pi} - \frac{1}{8} + \mathcal{O} \lk \frac{1}{T} \right) \right) 
+ \il_{C_T} \frac{\zeta' (s)}{\zeta (s)} d s \, .
\ee
\be
\boxed{\rule[-.3\baselineskip]{0pt}{9mm}
N (T) = 1 - \frac{T}{2 \pi} \log \pi + \frac{T}{2 \pi} \log \frac{T}{2} - 
\frac{T}{2 \pi} - \frac{1}{8} + \mathcal{O} \lk 
\frac{1}{T} \right) + \frac{1}{2 \pi i} \il_{C_T} \frac{\zeta' (s)}{\zeta (s)} \, . }
\ee
The last term can be split up into two parts, the results of which are given without further detailed calculations:
\be
\il^{2 + i T}_{2 - i T} \frac{\zeta' (s)}{\zeta (s)} d s = \mathcal{O} (1) \, , \quad \quad \mbox{for} \quad T \geq 3
\ee
and using 
\begin{align}
 \il^{2 - i T}_{\frac{1}{2} - i T} \frac{\zeta' (s)}{\zeta (s)} d s &= \il^2_{1/2} \frac{\zeta' (\sigma - i T)}{\zeta (\sigma - i T)} d \sigma = 
 \overline{\il^2_{1/2} \frac{\zeta' (\sigma - i T)}{\zeta (\sigma + i T)} } d s \nonumber\\
 &= \overline{\il^{2 + i T}_{\frac{1}{2} + i T} \frac{\zeta' (s)}{\zeta (s)} d s} \nonumber\\
 > \quad \quad & \frac{1}{2 \pi i} \Big( \il^{2 - i T}_{\frac{1}{2} - i T} \frac{\zeta' (s)}{\zeta (s)} d s + 
 \il^{\frac{1}{2} + i T}_{2 + i T} 
 \frac{\zeta' (s)}{\zeta (s)} d s \Big) =  \frac{1}{\pi} \,\mathrm{Im} \lk \il^{\frac{1}{2} + i T}_{2 + i T} 
 \frac{\zeta' (s)}{\zeta (s)} d s \right) \, .
\end{align}
So far we have found
\be
N (T) = \frac{T}{2 \pi} \log \frac{T}{2 \pi} - \frac{T}{2 \pi} + \frac{7}{8} + \mathcal{O} \lk \frac{1}{T} \right) + 
\frac{1}{\pi} Im \lk \il^{\frac{1}{2} + i T}_{2 + i T} \frac{\zeta' (s)}{\zeta (s)} d s \right) \, .
\ee
Using
\begin{align}
 \il^{\frac{1}{2} + i T}_{2 + i T} \frac{\zeta' (s)}{\zeta (s)} d s &= \log \zeta \lk \frac{1}{2} + i T \right) - 
 \log \zeta (2 + i T) \nonumber\\
 > \,\mathrm{Im} \lk  \il^{\frac{1}{2} + i T}_{2 + i T} \frac{\zeta' (s)}{\zeta (s)} d s \right) &= arg 
 \lk \zeta \lk \frac{1}{2} + i T \right) \right) - arg \lk \zeta (2 + i T) \right) \, .
\end{align}
The modulus of the last expression can be shown to be $\mathcal{O} (\log T)$.

Hence our final result for the number of zeros in the critical strip with $0<T$ is given by
\be
\boxed{N (T) = \frac{T}{2 \pi} \lk \log \frac{T}{2 \pi} - 1 \right) + \mathcal{O} (\log T) \, . } 
\ee
As mentioned above, this formula was given by Riemann in 1859, but only proved by von Mangoldt in 1905.

By the way, we can also approximate $\,\mathrm{Im}\log \Gamma (1/4 + it/2)$
and so obtain
\begin{align}
 \,\mathrm{Im} \left\{ \log \Gamma \lk \frac{1}{4} + \frac{it}{2} \right) \right\} &=
 \frac{t}{2} \log \lk \frac{t}{2} \right) - \frac{t}{2}  - \frac{\pi}{8} - \frac{t}{2} \log \pi + \mathcal{O} (t^{- 1}) \nonumber\\
 \mbox{i.e.} \quad \vartheta (t) &= \frac{t}{2} \log \lk \frac{t}{2 \pi} \right) - \frac{t}{2} - 
 \frac{\pi}{8} + \mathcal{O} (t^{- 1}) \, .
\end{align}
This brings us to the useful result
\be
N (T) = \frac{1}{\pi} \vartheta (T)  + 1 + \frac{1}{\pi} arg \zeta \lk \frac{1}{2} + i T \right)  \, , 
\ee
with
\be
\frac{1}{\pi} arg \zeta \lk \frac{1}{2} + i T \right) = \mathcal{O} (\log T) \quad \mbox{for} \quad T \to \infty \, .
\ee
So we can conclude for the number of zeros of $\zeta$ in the critical strip:
\begin{align}
 1. & N (T) \stackrel{T \to \infty}{\longrightarrow} \infty \nonumber\\
 2. & N (T) \sim \frac{T}{2 \pi} \log T  \, .
\end{align}
This follows from
\be
N (T) = \frac{T}{2 \pi} \log \frac{T}{2 \pi} + \mathcal{O} (\log T) \, ,
\ee
which when divided by $T/2 \pi \log T$, leads to
\be
\frac{N (T)}{\frac{T}{2 \pi} \log T} = \frac{\log T - \log 2 \pi}{\log T} + \frac{C}{T/2 \pi} 
\,\,\underset{T \to \infty}{\longrightarrow} 1 \, .
\ee
This result should be compared with the prime number theorem (Gau\ss{} 1796, when he was 15 years old)
\be
\pi (x) \sim \frac{x}{\log x} \quad  \mbox{or} \quad \lim\limits_{x \to \infty} \lk \frac{\pi (x)}{\frac{x}{\log x}} \right) = 1 \, .
\ee
Von Koch proved in 1901: If the Riemann hypothesis $\lk Re (s) = \frac{1}{2} \right)$ is true, then
\be
\pi (x) = Li (x) + \mathcal{O} \lk \sqrt{x} \log x \right) \, , 
\ee
i.e., the error in the claim $\pi (x) \sim Li (x)$ is of the order $\sqrt{x} \log x$. 
\bi

\no
\section{Riemann's Zeta Function Regularization}
\bi

\no
In this final section, we want to introduce the concept of the zeta function in connection with regularizing certain problems in quantum 
physics where infinities occur. For this reason, we consider an operator $A$
with positive, real discrete eigenvalues $\{a_n \}$, i.e., $ A f_ n(x) = a_nf(x)$ and one defines its associated zeta function by
\be
\zeta_A (s) = \sli_n a^{- s}_n = \sli_n e^{- s \ln a_n} \, ,
\ee
where $n$ runs over all eigenvalues. If one chooses for $A$
the Hamilton operator of the harmonic oscillator, for example, one gets (apart from the zero-point energy) exactly the Riemann zeta function. 
By formal differentiation now follows:
\be
\zeta'_A (0) = \left. - \sli_n \ln a_n e^{- s \ln a_n} \right|_{s = 0} = - \ln \lk \pli_n a_n \right) \, .
\ee
This suggests the definition
\be
\det A = \exp \left[ - \zeta'_A (0) \right] \, ,
\ee
which we shall exclusively be using in the following. The advantage of this method is that $\zeta'_A (0)$ 
is not singular for many operators of physical interest. As an example of the many applications to relativistic as well as 
non-relativistic problems in quantum field theory, we will choose the Casimir effect. 

This effect is a non-classical electromagnetic, attractive or repulsive force which occurs between electrically neutral conductors 
in a vacuum. The  size of this force was first calculated by Casimir for the case of ideal conducting, infinitely extended, 
parallel plates; his result was a force
\be
F = - \frac{\pi^2}{240} \cdot \frac{\hbar c}{a^4} \, ,
\ee
where $a$ is the distance between the plates and the negative sign indicates that the plates attract each other. 
This force apparently depends only on the fundamental constants $\hbar$ and $c$
apart from the distance between the plates; not, however, on the coupling constant $\alpha$
between the Maxwell and the matter field. Its quantum mechanical character is revealed by the fact that $F$
vanishes in the classical limit
$\hbar \to 0$.
\bi

\no
Casimir's derivation of $F$ was based on the concept of a quantum electrodynamic (particle) vacuum representing the zero-point 
oscillations of an infinite number of harmonic oscillators. As a result, one gets the total vacuum energy by summation over the 
zero-point energies $1/2 \hbar \omega_{\vk}$  of all allowed modes
with wave number vector $\vk$ and polarization $\sigma$,
\be
E = \sli_{\vk, \sigma} \frac{1}{2} \hbar \omega_{\vk} \, .
\ee
If we evaluate this equation for the case of two plane parallel plates at distance $a$ from each other, one does get a divergent 
total energy $E(a)$, but the energy difference $E(a)-E(a+ \delta a)$ is finite ($\delta a$ = infinitesimal change in the plate distance), 
leading also to a finite force per unit area,
\be
F = - \frac{\partial E (a)}{\partial a} \, .
\ee
To calculate this energy difference or force, a UV-cut-off is usually introduced, i.e., the energy $E$ is replaced by
\be
\sli_{{\vk}, \sigma} \frac{1}{2} \hbar \omega_{\vk} e^{- \frac{b}{\pi c} \omega_{\vk}}
\ee
and, in the end result, the limit $b \to 0$ is considered.
\bi

\no
This derivation of $F$, however, can give the impression that the appearance of the Casimir force is linked to the existence of the 
zero-point fluctuations of the quantized electromagnetic field.
\bi

\no
In order to avoid the divergent vacuum energy problem, in the following, we shall consider the problem according to Hawking from the 
viewpoint of path integral quantization and zeta-function regularization. Here, it is again unnecessary to refer to the vacuum oscillation. 
For reasons of simplicity, we wish to consider the Casimir effect only for a real, scalar field theory which is defined by $(\hbar  = c = 1!)$
\be
\cL (\phi) = - \frac{1}{2} \partial_\mu \phi \partial^\mu \phi - \frac{1}{2} m^2 \phi^2 - V (\phi) \, ,
\ee
with the arbitrary potential $V$.
\bi

\no
First, we couple the field $\phi$ to an external source $J$,
\be
\cL (\phi) \to \cL (\phi) + J \phi \, .
\ee
We can then write the vacuum amplitude $ \langle 0_+ | 0_- \rangle^J$  or the action $W[J]$ in the form
\be
\langle 0_+ | 0_- \rangle^J = e^{i W [J]} = \int [d \phi] e^{i \int d^4  x \{ \cL (\phi) + J \phi \}} \, ,
\ee
where we guarantee the convergence of the path integral by the substitution 
$m^2 \to m^2  -i \epsilon \, , \epsilon  > 0$. We have assumed that  $|0_- \rangle$ or $|0_+ \rangle$ 
describes a vacuum which is not ``disturbed'' by the presence of certain geometries, i.e., the path integral is, 
without restriction by boundary conditions, to be taken over all fields $\phi$. 
This changes as soon as we introduce two plates into the vacuum, for example, perpendicular to the $z$ axis 
(points of intersection: $z = 0$ and $z = a$) 
and require that only those fields should contribute to the path integral which would vanish on the plate surface, i.e., 
for which it holds that
\be
\phi (x_0, x_1, x_2, 0) = \phi (x_0, x_1, x_2, a) = 0
\ee
for arbitrary $(x_0, x_1, x_2)$. We now get
\begin{align}
 \langle 0_+ | 0_- \rangle^J_a & = e^{i W (a, [J])} \nonumber\\
 &= \il_{\cF_a} [d \phi] \exp \left[ i \int d^4 x \left\{ - \frac{1}{2} \partial_\mu \phi \partial^\mu \phi - 
 \frac{1}{2} (m^2 - i \epsilon) \phi^2 - V (\phi) - J \phi \right\} \right] \, ,
\end{align}
where $\int_{\cF_a}$ suggests that the path integral is only to be taken over the restricted space of functions $\cF_a$
defined by the boundary conditions. With this, we have represented the vacuum amplitude or the action for the most general case as a 
function of the geometric parameter $a$ and as a functional of the external source $J$. 
In order to approach the conditions of the QED Casimir effect, we now choose $J = 0$
as well as a free $(V = 0)$, massless $(m = 0)$ field $\phi$. Following a partial integration:
\be
\langle 0_+ | 0_- \rangle_a = e^{i W (a)} = \il_{\cF_a} [d \phi] e^{- \frac{i}{2} \int d^4 x \phi \{ - \partial^2 - i \epsilon \} \phi} \, .
\ee
The Gauss integral gives 
\begin{align}
\langle 0_+ | 0_- \rangle_a  = 
e^{i W (a)} &= \il_{\cF_a} [d \phi] 
e^{- \frac{1}{2} \int d^3 x d \tau \phi \{ - \Box {}_E \} \phi} \, .
\end{align}
Here, $N$ is a (divergent) constant which we shall set $= 1$, 
since it only contributes a non-physical additive constant to $W(a)$. By writing $\Box_E/\cF_a$, 
we mean that only eigenvalues with eigenfunctions in $\cF_a$ can be used to evaluate the determinant. Furthermore (in keeping with the 
$i \epsilon$ requirement), a Wick rotation $t \to i \tau$ was made, i.e., $\Box_E = \partial^2_\tau + \Delta$.
\bi

\no
From the original definition of the determinant, it follows that
\begin{align}
  \langle 0_+ | 0_- \rangle_a  = e^{i W (a)} &=  \left[ \exp \left\{ - \zeta'_{- \Box_E / \cF_a} (0) \right\} \right]^{- \frac{1}{2}} \nonumber\\
  &= \exp \left[ \frac{1}{2} \zeta'_{- \Box_E / \cF_a} (0) \right] \, .
\end{align}
The operator $ - \Box_E/\cF_a$ has the spectrum
 \be
 \left\{ k^2_0 + k^2_1 + k^2_2 + \lk \frac{\pi n}{a} \right)^2 | k_0, k_1, k_2 \in \RR, n \in \NN \right\}
\ee
and thus, the zeta function
\be
\zeta_{- \Box_E / \cF_a} (s) = 2 \frac{A}{(2 \pi)^2} \frac{T_E}{2 \pi} \int \il^\infty_{- \infty} \int d k_0 d k_1 d k_2 
\sli^\infty_{n = 1}
\left[ k^2_0 + k^2_1 + k^2_2 + \lk \frac{n \pi}{a} \right)^2 \right]^{- s} \, .
\ee
Here, the factor $2$ makes allowance for the two polarization possibilities of the photon, which, in our simple model, have no analogue. 
Furthermore, $AT_E$ is a normalization volume in three-dimensional $(0, 1, 2)$
space, where the Euclidean time $T_E$ is linked to a (Minkowski) normalization time interval $T$ by $T_E = iT$. 
Dropping the term independent of a $(n = 0)$ in the last equation simply leads to the subtraction of an (infinite) constant of $W(a)$.
\bi

\no
Further evaluation of $\zeta_{- \Box_E / \cF_a} (s)$ now takes on the form
\begin{align}
 \zeta_{- \Box_E / \cF_a} (s) &= 2 A T_E \frac{4 \pi}{(2 \pi)^3} \sli^\infty_{n = 1} \il^\infty_0 d k k^2 \left[ k^2 + \lk \frac{n \pi}{a} 
 \right)^2 \right]^{- s} \nonumber\\
 &= \frac{8 \pi}{(2 \pi)^3} A T_E \lk \frac{\pi}{a} \right)^{3 - 2s} \sli^\infty_{n = 1} n^{3 - 2s} \frac{1}{2}
 \frac{\Gamma \lk \frac{3}{2} \right) \Gamma \lk s - \frac{3}{2} \right)}{\Gamma (s)} \nonumber\\
 &= \frac{4 \pi}{(2 \pi)^3} A T_E \lk \frac{\pi}{a} \right)^{3 - 2s} \zeta (2 s - 3) \frac{\Gamma 
 \lk \frac{3}{2} \right) \Gamma \lk s - \frac{3}{2} \right)}{\Gamma (s)}  \, .
\end{align}
The derivative is
\begin{align}
 \zeta'_{- \Box_E / \cF_a} (0) &= \left. \frac{4 \pi}{(2 \pi)^3} A T_E \lk \frac{\pi}{a} \right)^3 \zeta (- 3) 
 \Gamma \lk \frac{3}{2} \right) \Gamma \lk - \frac{3}{2} \right) \frac{d}{ds} \frac{1}{\Gamma (s)} \right|_{s = 0} \nonumber\\
 &= \frac{\pi^2}{360 a^3} A T_E \, .
\end{align}
Finally we get
\be
\langle 0_+ | 0_- \rangle = e^{i W (a)} = e^{- \epsilon (a) T_E} = e^{- i \epsilon (a) T} \, ,
\ee
with
\be
\epsilon (a) = - \frac{\pi^2}{720 a^3} A \, .
\ee
The appearance of the phase factor $e^{-i \epsilon (a)T}$ 
in the vacuum amplitude allows us to identify $\epsilon (a)$ 
as the vacuum energy displacement and to write, for the force per surface unit,
\be
F = - \frac{1}{A} \frac{\partial \epsilon}{\partial a} \, ,
\ee
which leads to 
\be
F = - \frac{\pi^2}{240} \cdot \frac{1}{a^4} 
\ee
or, after putting  $\hbar$ and $c$ back in:
\be
F = - \frac{\pi^2}{240} \cdot \frac{\hbar c}{a^4} \, .
\ee
This is precisely Casimir's result which we have now completely derived with the aid of 
Riemann's zeta-function regularization, 
which completely eliminated the divergent zero-point energy. The same procedure finds application in QED and QCD, and can be looked up 
in the list of references (i.e., in \cite{9,10,11}).

\section*{Acknowledgment}

I wish to express my sincere gratitude to the librarians at the ``Handschriftenabteilung''
(Department of Handwritten Documents) at  G\"ottingen University for giving me access to Riemann's original handwritten manuscripts, 
in particular to the originals concerning prime numbers.

\appendix

\part*{Supplements}

The Riemann $\zeta$ function can be extended meromorphically into the region $\lbrace s:\Re(s)>0\rbrace$ in and on the right of the critical strip
$\lbrace s: 0\leq \Re(s)<1\rbrace$. This is a sufficient region of meromorphic continuation for many applications in analytic number theory.
The zeroes of the $\zeta$ function in the critical strip are known as the non-trivial zeroes of $\zeta$.

It is remarkable that $\zeta$ obeys a functional equation establishing a symmetry across the critical line $\lbrace s:\Re(s)=\tfrac{1}{2}\rbrace$
rather than the real axis. One consequence of this symmetry is that the $\zeta$ function may be extended meromorphically to the entire complex plane with a simple pole at $s=1$
and no other poles. For all $\mathfrak{C}\setminus \Re(s)=1$ including the strip we have the functional equation:
\begin{equation}
\zeta(s) = 2^s\pi^{s-1}\sin\left(\frac{s\pi}{2}\right) \Gamma(1-s)\zeta(1-s),\quad \Re(s)<0 \label{eq:RiemannZetaFunctionalEquation}
\end{equation}
or, equivalently, the identity between meromorphic functions $\zeta(s)$:
\begin{equation}
\zeta(1-s) = \frac{2}{(2\pi)^s}\cos\left(\frac{s\pi}{2}\right)\Gamma(s)\zeta(s).\label{eq:RiemannZetaMeromorphicFunctionalEquation}
\end{equation}

The analytical continuation given here allows one to connect $\zeta(s)$ for positive values of $\Re(s)$ with the same for negative values, for instance:
\begin{align}
\zeta(-1) = 2^{-1}\pi^{-2}(-1)\Gamma(2)\zeta(2) &= \frac{1}{2}\cdot\frac{1}{\pi^2}\cdot(-1)\cdot 1\cdot\frac{\pi^2}{6} = -\frac{1}{12}\; , \label{eq:SymmetryExample}
\intertext{i.e.,}
\zeta_R(-1) &= -\frac{1}{12}\quad , \label{eq:RiemannZetaAtNegativeUnity}
\end{align}
where the subscript $R$ is added to distinguish Riemann's $\zeta$ from Euler's $\zeta$, of which it is an extension, i.e.,
\begin{align*}
\zeta(x) &= \sum_{n=1}^\infty \frac{1}{n^x} = \prod_{p \text{ prime}}\frac{1}{1-p^{-x}}\quad \text{converging for } x>1\\
\frac{1}{1^x}+\frac{1}{2^x}+\frac{1}{3^x}+\dots &= \prod_{p \text{ prime}}\frac{p^x}{p^x-1} = \left(\frac{2^x}{2^x-1}\right)\left(\frac{3^x}{3^x-1}\right)\left(\frac{5^x}{5^x-1}\right)\dots
\end{align*}

\begin{figure}[!hbt]
\centering
\includegraphics[width=0.7\linewidth]{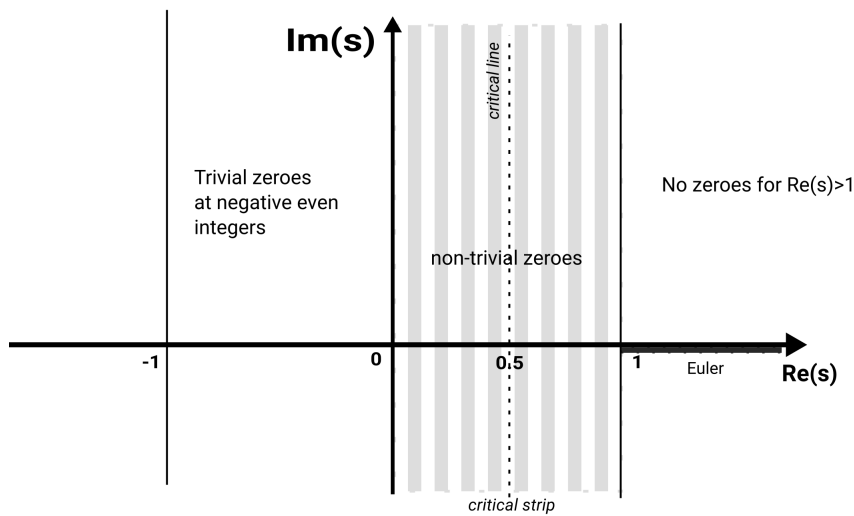}
\caption{The different domains of definition of Riemann's $\zeta$ function of \eqref{eq:RiemannZetaInDifferentDomains} }
\label{fig:RiemannZetaDomains}
\end{figure}


When we extend this function into the whole complex $s$ plane, then Riemann's $\zeta$ function comes in three different representations:
\begin{equation}
\zeta(s) = \begin{cases}
\sum_{n=1}^\infty \frac{1}{n^s} = \prod_{p \text{ prime}}\frac{p^s}{p^s-1}\; , & \Re(s) >1\\
(1-2^{1-s}) \sum_{n=1}^\infty \frac{(-1)^{n+1}}{n^s}\;, & 0<\Re(s)<1\\
2^s\pi^{s-1}\sin\left(\frac{s\pi}{2}\right) \Gamma(1-s)\zeta(1-s)\;, & \Re(s)<0
\end{cases}
\label{eq:RiemannZetaInDifferentDomains}
\end{equation}

Where is $\zeta(s)$ equal to zero?
\begin{enumerate}
\item No zeroes for $\Re(s)>1$ since here $\zeta(s)>0$.
\item Non-trivial zeroes in the strip $0<\Re(s)<1$, symmetric around $\Re(s)=\tfrac{1}{2}$.
\item Trivial zeroes for $s=-2,-4,\dots$, thus for $\Re(s)<0$ .
\end{enumerate}
There is a pole at $s=1$.

\section*{The origins of the functional equation for Dirichlet's $\eta$ function}
Euler in his ``Remarques sur un beau rapport entre les series des puissances tant directes que reciproches''\footnote{Remarks on a beautiful relation between direct as well as reciprocal power series.}  writes the following functional equations
\begin{align*}
\frac{1-2^{n-1}+3^{n-1}-4^{n-1}+5^{n-1}-6^{n-1}+\dots}{
1-2^{-n}+3^{-n}-4^{-n}+5^{-n}-6^{-n}+\dots}
&= -\frac{1\cdot 2\cdot 3 \cdot \dots(n-1)(2^n-1)}{(2^{n-1}-1)\pi^n}\cos\left(\frac{n\pi}{2}\right)\\
\frac{1-3^{n-1}+5^{n-1}-7^{n-1}+\dots }{ 1-3^{-n}+5^{-n}-7^{-n}+\dots }
&= \frac{1\cdot 2\cdot 3 \cdot \dots (n-1)(2^n)}{\pi^n}\sin\left(\frac{n\pi}{2}\right)\quad .
\end{align*}
Then he finishes his work by proving that the above statements hold true for positive and negative whole numbers as well as for fractional values of $n$.

Nowadays we write with $s\in\mathbb{C}$:
\begin{equation}
\eta(1-s) = -\frac{(2^s-1)}{\pi^s(2^{s-1}-1)}\cos\left(\frac{\pi s}{2}\right)\Gamma(s)\eta(s)
\end{equation}
which is the functional equation of Dirichlet's $\eta$ function.

Hardy gave a proof for the case when $s$ is replaced by $s+1$ in the last equation:
\begin{equation}
\eta(-s) = 2\frac{(1-2^{-s-1}}{1-2^{-s} }\pi^{-s-1} s \sin\left(\frac{\pi s}{2}\right)\Gamma(s)\eta(s+1)\quad .
\end{equation}
From the relation $\eta(s) = \left(1-2^{1-s}\right)\zeta(s)$ one can show that $\eta$ has zeroes at the points 
$s_k=1+\frac{2\pi i k }{\ln 2}$ for all $k\in\mathbb{Z}\setminus\lbrace 0\rbrace$, e.g., $s_1 = 1+9.0647i$. For $k=0$ one finds instead $\eta(1) = \ln 1 = 0.69315$.
Remember that $\zeta(1) = \infty$.

When we write
\[ \zeta(s) = \frac{\eta(s)}{1-2^{1-s}} \]
we realize that $\eta(s)$ as well as $(1-2^{1-s})$ have the same zeroes $s_k$ with $k=1,2,3,\dots$.
$\eta(s)$ is also zero at the points where $\zeta(s)$ is zero. These are the trivial zeroes $s=-2,-4,-6,\dots$ such that
\[ \eta(-2) = \eta(-4)=\eta(-6) = \dots = 0\quad .
\]
Finally, $\eta$, like $\zeta$, possesses the non-trivial zeroes within the critical strip $\lbrace s\in\mathbb{C} \vert 0<\Re(s)<1\rbrace$.
The celebrated unproven Riemann hypothesis claims that all non-trivial zeros of $\zeta$ are located on the axis $\Re(s) = \tfrac{1}{2}$.

\underline{$\zeta(s)$ is  a meromorphic function}.
Later we will meet Riemann's $\xi$ function, $\xi(s)=\tfrac{1}{2}s(s-1)\pi^{-\tfrac{s}{2}}\Gamma\left(\frac{s}{2}\right)\zeta(s)$. $\xi(s)$ is an \underline{entire function}, it has non-trivial
zeroes, however no trivial zeroes and no poles. Also: $\xi(s)=\xi(1-s)$.

\begin{table}[!tbp]

\centering
 \subfloat[A few values of $\zeta(s)$]{
$\begin{array}{|lr|}\hline
s & \zeta(s)\\\hline\hline
-2\mathbb{N} & 0\\
-\mathbb{N}  &\frac{-B_{n+1}}{n+1} \\
-7 & \frac{1}{240}\\
-5 & \frac{-1}{252}\\
-3 &  \frac{1}{20}\\
-1 & \frac{-1}{12}\\
0 &  -\frac{1}{2}\\
\frac{1}{2} & -1.46035450\\
1 & \infty \\
\frac{3}{2} & 2.6123753486\\
2  & \frac{\pi^2}{2}\approx 1.6449340 \text{(Euler,Basel)} \\
\frac{5}{2} & 1.3414872572\\
3 & 1.2020569\\
\frac{7}{2} & 1.1267338673 \\
4 & \frac{\pi^4}{90}\approx 1.082323233\\\hline
\end{array}
$
}\;
\subfloat[Function Properties]{
\begin{tabular}{|lccc|}\hline
Function & Non-trivial zeroes & Trivial zeroes & Poles\\\hline
$\zeta(s)$ & Yes & $-2,-4,-6,\dots$ & $1$\\
$\zeta(1-s)$ & Yes & $3,5,\dots$ & $0$ \\
$\sin\frac{\pi s}{2}$ & No & $2\mathbb{N}$ & No\\
$\cos\frac{\pi s}{2}$ & No & $2\mathbb{N}+1$ & No\\
$\sin\pi s$ & No & $\mathbb{N}$ & No\\
$\Gamma(s)$ & No & No & $0,-1,-2,\dots$\\
$\Gamma\left(\frac{s}{2}\right)$ & No & No & $0,-2,-4,\dots$\\
$\Gamma(1-s)$ & No & No & $1,2,3,\dots$\\
$\Gamma\left(\frac{1-s}{2}\right)$ & No & No & $1,3,5,\dots$\\
$\xi(s)$ & Yes & No & No\\\hline
\end{tabular}
}
\caption{Properties and special values of the Riemann $\zeta$ function.}
\label{tbl:PropertiesOfRelevantFunctions}
\end{table}

\begin{figure}[tp]
\centering
\includegraphics[width=0.8\linewidth]{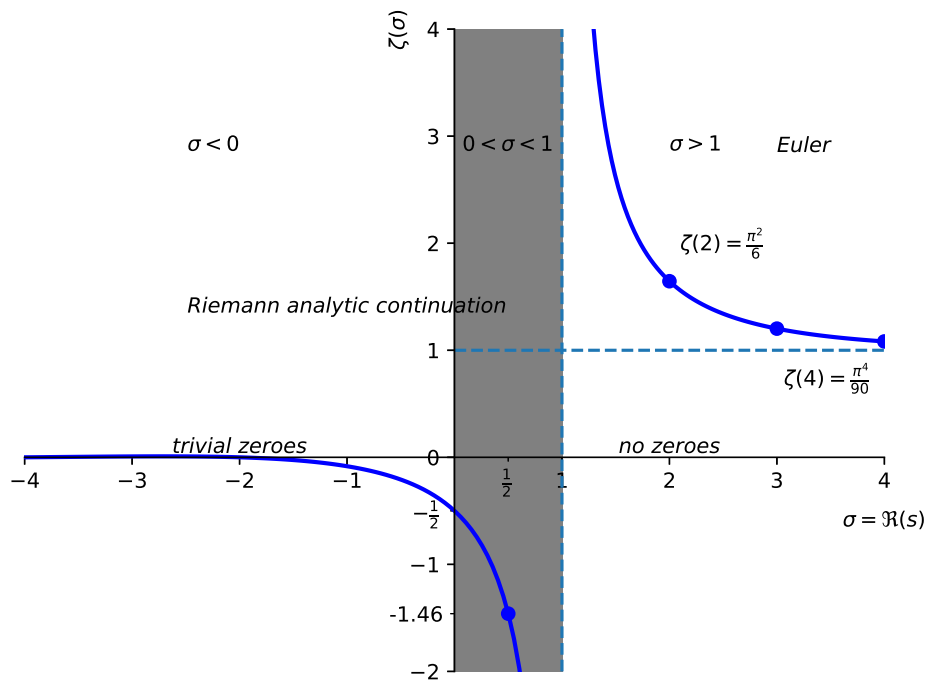}
\caption{The behaviour of Riemann's $\zeta$-function for real arguments.}
\label{fig:ZetaRealArgumentBehaviour}
\end{figure}

The tables (\ref{tbl:PropertiesOfRelevantFunctions}) indicate that the $\Gamma$ function and trigonometric factors in the functional equation (\eqref{eq:RiemannZetaFunctionalEquation}, \eqref{eq:RiemannZetaMeromorphicFunctionalEquation}, resp.)
are tied to the trivial zeros and poles of the $\zeta$ function, but have no direct bearing on the distribution of the non-trivial zeroes, which is the most important feature of the $\zeta$ function
for the purposes of analytic number theory, beyond the fact that they are symmetric about the real axis and the critical line $x=\tfrac{1}{2}$.
Exponential functions such as $2^{s-1}$ or $\pi^{-s}$ have neither zeroes nor poles. In particular the Riemann hypothesis is not going to be resolved just from further analysis of 
the $\Gamma$ function.

\underline{Remarkable historical fact}:
Euler, in 1749 (~110 years before Riemann!) discovered that the following series is convergent:
\begin{equation}
\phi(s) = \sum_{n=1}^\infty \frac{(-1)^{n+1}}{n^s}\label{eq:EulerPhi}
\end{equation}
This is also referred to as Dirichlet's $\eta$ function. It is related to $\zeta$ by
\begin{equation}
\phi(s) = (1-2^{1-s})\zeta(s)\label{eq:RelationRiemannZetaEulerPhi}
\end{equation}

Within the critical strip $0<s<1$ we have:
\begin{align}
\zeta(s) &= \frac{2^{s-1}}{2^{s-1}-1}\phi(s) = \frac{1}{1-2^{1-s}}\phi(s) \nonumber \\
&= \frac{1}{1-2^{1-s}}\sum_{n=1}^\infty \frac{(-1)^{n+1}}{n^s},\quad \Re(s) >0, 1-2^{1-s}\neq 0\;.
\end{align}

From Euler we have
\begin{equation}
\frac{\phi(1-n)}{\phi(n)} = \frac{-(n-1)!(2^n-1)}{(2^{n-1}-1)\pi^n}\cos\left(\frac{n\pi}{2}\right)\;,
\end{equation}
and he furthermore says: ``I shall hazard the following conjecture:
\begin{equation}
\frac{\phi(1-s)}{\phi(s)} = -\frac{\Gamma(s)(2^s-1)\cos\left(\frac{\pi s}{2}\right)}{(2^{s-1}-1)\pi^s}\label{eq:EulerConjectureOnPhi}
\end{equation}
is true for all $s$''.
We know that $(\eta(s) =)\phi(s) = (1-2^{1-s})\zeta(s)$, which leads at once from \eqref{eq:EulerConjectureOnPhi} to 
\begin{equation}
\zeta(1-s) = \frac{2}{(2\pi)^s}\Gamma(s)\zeta(s)\cos\left(\frac{\pi s}{2}\right),\quad \forall s\in\mathbb{C}\setminus 1
\end{equation}
and this is the famous functional equation which was proven by Riemann in 1859 (but it was conjectured by Euler in 1749!). It is probably correct to assume that Riemann was very familiar
with Euler's contribution.


With the alternating Dirichlet series at hand we can already make an important statement regarding the zeroes of the $\zeta$ function within the critical strip $0<\Re(s)=\sigma <1$,
which is important for the Riemann hypothesis, which claims that all non-trivial zeroes of $\zeta$ lie on the line with $\Re(s) = \tfrac{1}{2}$.

To show this we start with
\begin{equation}
\zeta(s) = \sum_{n=1}^\infty \frac{1}{n^s},\quad s:=\sigma+ i t\label{eq:RiemannZetaSeriesRepresentation}
\end{equation}
which is convergent for $\Re(s) > 1$, is a meromorphic function and has a pole at $s=1$.
Next let 
\begin{equation}
n^s = n^{\sigma + it} = n^\sigma n^{it} = n^\sigma e^{it\ln n} 
= \vert n\vert^\sigma \left( \cos(t\ln n) + i \sin(t\ln n)\right)\label{eq:AnalyticContinuationOfTheArgument}
\end{equation}
from which immediately follows
\begin{align}
\zeta(s) &= \Re(\zeta(s))+i\Im(\zeta(s)) = \sum_{n=1}^\infty \frac{1}{n^\sigma}\left[ \cos( t\ln n) - i\sin(t\ln n)\right]\\
\Rightarrow \Re(\zeta(s)) &= \sum_{n=1}^\infty n^{-\sigma}\cos(t\ln n)\label{eq:RiemannZetaRealPart}\\
\Im(\zeta(s)) &= \sum_{n=1}^\infty n^{-\sigma} \sin(t \ln n)\label{eq:RiemannZetaImagPart}
\end{align}
which are convergent for $\sigma >1, t\in\mathbb{R}$. Next consider the Euler's $\phi$ function as given in \eqref{eq:EulerPhi}, which is also known as
Dirichlet's $\eta$ function. An extension of the domain of $\zeta$ into the region of $0 < \sigma < 1$, i.e., into the critical strip, is obtained by
rewriting \eqref{eq:RelationRiemannZetaEulerPhi} as
\begin{equation}
\zeta(s) = \frac{1}{1-2^{1-s}}\eta(s)\; .
\end{equation}
Note that only the critical strip is of importance for the Riemann hypothesis. Note further that $\eta$ is convergent for $\sigma=\Re(s) > 0$ and that the following
alternating harmonic series,
\begin{equation}
\eta(1) = 1-\frac{1}{2}+\frac{1}{3}-\frac{1}{4}+\dots = \ln 2 \approx 0.69315\; ,\label{eq:DirichletEtaAtUnity}
\end{equation}
is obtained from
\begin{equation}
\ln(x+1) = x -\frac{1}{2}x^2+\frac{1}{3}x^3-\dots \quad -1<x\leq 1\;,
\end{equation}
where $x$ is assumed to be real.
One may rewrite Dirichlet's $\eta$ function in the following way:
\begin{equation}
\eta(s) = \sum_{n=1}^\infty\left(\frac{1}{(2n-1)^s} - \frac{1}{(2n)^s}\right)\quad .\label{eq:DirichletEtaSeriesExtension}
\end{equation}
From which one then obtains in a simple way (c.f. \eqref{eq:RiemannZetaRealPart}, \eqref{eq:RiemannZetaImagPart}):
\begin{align}
\Re(\eta(s)) &= \sum_{n=1}^\infty\left[(2n-1)^{-\sigma}\cos(t\ln(2n-1)) - (2n)^{-\sigma}\cos(t\ln(2n))\right]\label{eq:DirichletEtaSeriesRealPart}\\
\Im(\eta(s)) &= \sum_{n=1}^\infty\left[(2n)^{-\sigma}\sin(t\ln(2n))-(2n-1)^{-\sigma}\sin(t\ln(2n-1))\right]\quad .\label{eq:DirichletEtaSeriesImagPart}
\end{align}
Using $\cos x - \sin x = \sqrt{2}\sin\left(x+\tfrac{3}{4}\pi\right)$ one then obtains
\begin{multline}
\Re(\eta(s))+\Im(\eta(s)) = \sqrt{2}\sum_{n=1}^\infty\left[(2n-1)^{-\sigma}\sin\left(t\ln(2n-1)+\frac{3}{4}\pi\right)\right.\\
\left.-(2n)^{-\sigma}\sin\left(t\ln(2n)+\frac{3}{4}\pi\right)\right] \neq 0\quad \forall\sigma\in(0,\frac{1}{2}),\forall t\quad ,
\end{multline}
i.e., $\eta$ possesses no roots on the left half of the critical strip, and because of the reflection formula \eqref{eq:RiemannZetaMeromorphicFunctionalEquation}
this holds true for the right half as well, i.e., they can only be on the critical line $\sigma=\tfrac{1}{2}$, which is the \textbf{Riemann hypothesis}.

\begin{figure}
\centering
\subfloat[The argument]{
\framebox{
\includegraphics[width=0.5\linewidth]{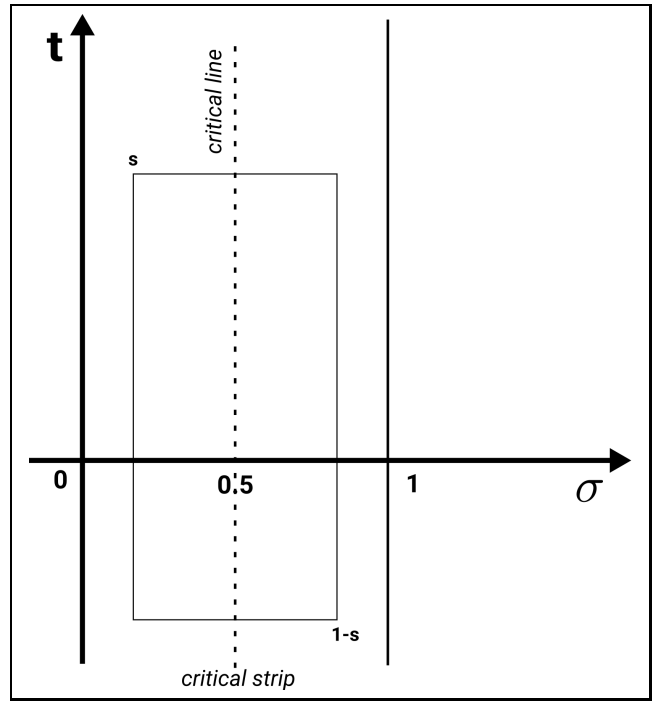}
}
}\;
\subfloat[The Function]{
\framebox{
\includegraphics[width=0.5\linewidth]{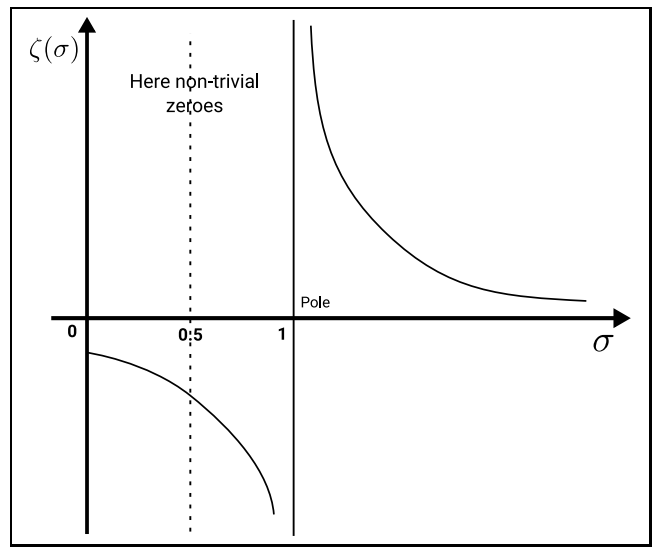}
\label{subfig:FunctionSketch}
}
}
\caption{A closer look at the behavior of $\zeta$. Referring to \ref{subfig:FunctionSketch} we have
$|\zeta(\tfrac{1}{2}-\sigma)|>|\zeta(\tfrac{1}{2}+\sigma)|\text{ or }|\zeta(\tfrac{1}{2}-\sigma)|>|\zeta(\tfrac{1}{2})|$.
No zeroes of $\zeta$ on the left half and right half of the critical strip, which is equivalent to Riemann's hypothesis.}
\end{figure}



\noindent \textbf{Theorem.} 
If $\Re(s) = \sigma > 0$ we have 
\begin{equation}
(1-2^{1-s})\zeta(s) = \eta(s) = \sum_{n=1}^\infty \frac{(-1)^{n-1}}{n^s}\;,\label{eq:ThmZetaFunctionRepresentation}
\end{equation}
which implies that $\zeta(s)<0$ if $s$ is real and $0<s<1$.

\begin{proof}
First assume that $\sigma > 1$ (Euler: $\Re(s)>1$). Then we have
\begin{align*}
(1-2^{1-s})\zeta(s) &= \sum_{n=1}^\infty \frac{1}{n^s}-2\sum_{n=1}^\infty \frac{1}{(2n)^s}\\
&=(1+2^{-s}+3^{-s}+\dots) - 2(2^{-s}+4^{-s}+6^{-s}+\dots)\\
&= 1-2^{-s}+3^{-s}-4^{-s}+\dots =\text{alternating $\zeta$ function},
\end{align*}
which proves \eqref{eq:ThmZetaFunctionRepresentation} for $\Re(s) = \sigma > 1$.
However, if $\sigma > 0$ the series on the right converges, thus \eqref{eq:ThmZetaFunctionRepresentation} also holds for $\sigma > 0$ by analytic continuation,
i.e., when $s$ is real then the sum in \eqref{eq:ThmZetaFunctionRepresentation} is an alternating series with a positive limit.

If $0<s<1$, then the factor $1-2^{1-s}$ becomes negative. Hence $\zeta(s)$ is also negative (has no zeroes!) in $0<s<1$.
\end{proof}

Note that $\eta(1)=\dots = \ln 2\approx 0.69315$ (c.f. \eqref{eq:DirichletEtaAtUnity}) while $\zeta(1)=\infty$, that is, $s=1$ is a pole of the meromorphic function $\zeta$.
Furthermore we have
\begin{equation}
\zeta(0) = -\frac{1}{2}\quad.\label{eq:ZetaValueAtOrigin}
\end{equation}
\begin{proof}
Starting with the functional equation
\begin{equation}
\Gamma\left(\frac{s}{2}\right)\pi^{-\frac{s}{2}}\zeta(s) = \Gamma\left(\frac{1-s}{2}\right)\pi^{-\frac{1-s}{2}}\zeta(1-s)
\end{equation}
solve for $\zeta(s)$ to obtain
\begin{align*}
\zeta(s) &= \pi^{\frac{s}{2}}\pi^{-\frac{1-s}{2}}\Gamma\left(\frac{1-s}{2}\right)\frac{\zeta(1-s)}{\Gamma\left(\frac{s}{2}\right)}\\
s\rightarrow 0\;:\; \zeta(0) &= \pi^{-\frac{1}{2}}\Gamma\left(\frac{1}{2}\right)\lim_{s\rightarrow 0}\frac{\zeta(1-s)}{\Gamma\left(\frac{s}{2}\right)}\;.
\end{align*}
Since the residues of $\zeta$ at $s=1$ and of $\Gamma$ at $s=0$ are both $1$, i.e.,
\begin{equation}
\zeta(s) = \frac{1}{s-1}+\dots,\quad \Gamma(s) = \frac{1}{s}+\dots\quad ,
\end{equation}
we have
\begin{equation}
\zeta(1-s) = -\frac{1}{s}+\dots,\quad \Gamma(\frac{s}{2}) = \frac{2}{s}+\dots
\end{equation}
and therefore
\begin{equation}
\lim_{s\rightarrow 0} \frac{\zeta(1-s)}{\Gamma(\frac{s}{2})} = \lim_{s\rightarrow 0}-\frac{\frac{1}{s}+\dots}{\frac{2}{s}+\dots} = -\frac{1}{2} 
\end{equation}
from which follows, using $\Gamma\left(\tfrac{1}{2}\right)$
\begin{equation}
\zeta(0) = \pi^{-\frac{1}{2}}\pi^{\frac{1}{2}}\left(-\frac{1}{2}\right) = -\frac{1}{2}\;\Longrightarrow\;  \zeta(0) = -\frac{1}{2}\;.
\end{equation}
\end{proof}

From the eqs.~(\ref{Gl: 98}), (\ref{Gl: 99}) we have
\[
t,x,\psi(x),\ln(x) \in\mathbb{R}\;.
\]
Therefore $\Im\xi(\tfrac{1}{2}+it)=0$, i.e., $\xi(\tfrac{1}{2}+it)\equiv \Xi(t) \in\mathbb{R}$
and thus
\begin{align*}
\Xi(t) = \xi(\frac{1}{2}+it) &=  -\frac{t^2+\frac{1}{4} }{ 2\left( \sqrt{\pi}\right)^{\frac{1}{2}+it} }
\Gamma\left( \frac{1}{4}+\frac{it}{2} \right)\zeta\left(\frac{1}{2}+it\right)\\
\xi(\frac{1}{2}) &= -\frac{1}{8\pi^{\frac{1}{4}}}\Gamma(\frac{1}{4})\zeta(\frac{1}{2})\approx 0.4971207781 =:a_0\\
\zeta(\frac{1}{2}) &\approx -1.4603545088\\
\Gamma(\frac{1}{4}) &= \sqrt{2\bar{\omega}2\pi} \approx 3.6256099082
\end{align*}
where in the last equation $\bar{\omega}$ is the so-called Gaussian lemniscate constant.

Some special values:
\begin{equation}
\xi(0) = \xi(1) =-\zeta(0) = \frac{1}{2}\quad .
\end{equation}

\begin{proof}
Using $\xi(s) = \frac{1}{2}s(s-1)\pi^{-\frac{s}{2}}\Gamma(\frac{s}{2})\zeta(s)$ as well as
$\Gamma(1+\frac{s}{2})=\frac{s}{2}\Gamma(\frac{s}{2})$ we obtain
\begin{equation}
\xi(s)\vert_{s=0} = (s-1)\pi^{-\frac{s}{2}}\Gamma(1+\frac{s}{2})\zeta(s)\vert_{s=0} \Leftrightarrow \xi(s) = -1\cdot 1\cdot \Gamma(1)\cdot \zeta(0) =\frac{1}{2}
\end{equation}
Thus 
\begin{equation*}
\xi(0) = \frac{1}{2}\quad .
\end{equation*}
In a similar manner, utilizing the reflection property $\xi(s)=\xi(1-s)$:
\begin{multline}
\xi(s) = (-s)\pi^{-\frac{1}{2}(1-s)}\Gamma(\frac{3}{2}-\frac{s}{2})\zeta(1-s)\\
\Rightarrow \xi(1) = -1 \cdot 1 \cdot \Gamma(1)\cdot \zeta(0) = \frac{1}{2}\\
\Longrightarrow \xi(1) = \frac{1}{2}
\end{multline}
\end{proof}

\section*{Riemann's Functional Equation}
\begin{equation}
\pi^{-\frac{s}{2}}\Gamma(\frac{s}{2})\zeta(s) = \pi^{-\frac{1-s}{2}} \Gamma(\frac{1-s}{2})\zeta(1-s) \, ,
\end{equation}
whose symmetry is obvious when $s\rightarrow 1-s$ is substituted into both sides of the equation.

\begin{proof}
Starting with Euler's $\Gamma$ function
\begin{equation}
\Gamma(s) = \int_0^\infty t^{s-1}e^{-t}dt\; .
\end{equation}
Using $s\rightarrow\tfrac{s}{2}$, the above results in
\begin{equation}
\Gamma(\frac{s}{2}) = \int_0^\infty t^{\frac{s}{2}-1}e^{-t}dt\;.
\end{equation}
Next one can use the substitution $t=\pi n^2x$ ($dt=\pi n^2 dx$) to obtain
\begin{align*}
\Gamma(\frac{s}{2}) &= \int_0^\infty (\pi n^2 x)^{\frac{s}{2}-1} e^{-\pi n^2 x}\pi n^2 dx\\
\pi^{-\frac{s}{2}}\Gamma(\frac{s}{2})\frac{1}{n^s} &= \int_0^\infty x^{\frac{s}{2}-1}e^{-\pi n^2 x} dx \, .
\intertext{Summation over $n$ yields}
\sum_{n=1}^\infty \pi^{-\frac{s}{2}}\Gamma(\frac{s}{2})\frac{1}{n^s} &= \sum_{n=1}^\infty\int_0^\infty x^{\frac{s}{2}-1}e^{-\pi n^2 x} dx\\
\pi^{-\frac{s}{2}}\Gamma(\frac{s}{2})\sum_{n=1}^\infty\frac{1}{n^s} &= \int_0^\infty x^{\frac{s}{2}-1} \sum_{n=1}^\infty e^{-\pi n^2 x} dx\\
\pi^{-\frac{s}{2}}\Gamma(\frac{s}{2})\zeta(s) &= \int_0^\infty x^{\frac{s}{2}-1}  \underbrace{ \sum_{n=1}^\infty e^{-\pi n^2 x} }_{\text{closely related to Jacobi $\vartheta$ func.}} dx\\
\vartheta (x) = \sum_{n\in\mathbb{Z}} e^{-\pi n^2 x} &= 1 + 2\sum_{n=1}^\infty e^{-\pi n^2 x} = 1+2\psi(x),\; x>0\; .\\
\Rightarrow \int_0^\infty x^{\frac{s}{2}-1}\sum_{n=1}^\infty e^{-\pi n^2 x} dx &= \int_0^\infty x^{\frac{s}{2}-1} \psi(x) dx\;.
\end{align*}
Split the integral on the r.h.s into two parts:
\begin{equation}
\int_0^\infty x^{\frac{s}{2}-1}\psi(x) dx = \int_1^\infty x^{\frac{s}{2}-1}\psi(x)dx+\int_0^1 x^{\frac{s}{2}-1}\psi(x)dx \, .
\end{equation}

Look at $\vartheta(x) = \frac{1}{\sqrt{x}}\vartheta(\tfrac{1}{x})$ or $2\psi(x)+1 = \tfrac{1}{\sqrt{x}}(1+\psi(\tfrac{1}{x}))$. The equations (\ref{Gl: 72})ff.~in the body of the paper are
\begin{align*}
\psi(x) &= \frac{1}{\sqrt{x}}\psi(\frac{1}{x}) -\frac{1}{2}+\frac{1}{2\sqrt{x}}\\
\int_0^1x^{\frac{s}{2}-1}\psi(x) dx &= \int_0^1 x^{\frac{1}{2}-1}\left(\frac{1}{\sqrt{x}}\psi(\frac{1}{x})+\frac{1}{2\sqrt{x}}-\frac{1}{2}\right)dx\\
&= \int_0^1\left(x^{\frac{s}{2}-\frac{3}{2}}\psi(\frac{1}{x})+\frac{1}{2}\left(x^{\frac{s}{2}-\frac{3}{2}}-x^{\frac{s}{2}-1}\right)\right) dx\\
&= \int_0^1 x^\frac{s-3}{2}\psi(\frac{1}{x})dx + \frac{1}{2}\left[\frac{1}{\frac{s}{2}-\frac{1}{2}}x^{\frac{s}{2}-\frac{1}{2}} -\frac{1}{ \frac{s}{2} }x^\frac{s}{2} \right]_0^1\\
&= \int_0^1 x^{\frac{s}{2}-\frac{3}{2}}\psi(\frac{1}{x})dx + \frac{1}{s(s-1)}\\
&\overset{(*)}{=} \int_\infty^1 \left( \frac{1}{y} \right)^{\frac{s}{2}-\frac{3}{2}} \psi(y) \left( -\frac{1}{y^2}\right)dy + \frac{1}{s(s-1)}\\
&\overset{y\rightarrow x}{=} \int_1^\infty \left( \frac{1}{x} \right)^{\frac{s}{2}-\frac{3}{2}} \psi(x) \frac{dx}{x^2} + \frac{1}{s(s-1)}\\
\Rightarrow \int_0^1 x^{ \frac{s}{2}-1 }\psi(x)dx &= \int_1^\infty x^{ -\frac{s}{2}-\frac{1}{2} }\psi(x) dx + \frac{1}{s(s-1)}\\
\int_0^\infty x^{\frac{s}{2}-1} \psi(x) dx &= \int_1^\infty x^{\frac{s}{2}-1}\psi(x)dx+\int_0^1x^{\frac{s}{2}-1}\psi(x)dx\\
&= \int_1^\infty x^{\frac{s}{2}-1}\psi(x)dx + \int_1^\infty x^{-\frac{s}{2}-\frac{1}{2}}\psi(x)dx + \frac{1}{s(s-1)}\\
&= \int_1^\infty\left(x^{\frac{s}{2}-1}+x^{-\frac{s}{2}-\frac{1}{2}}\right)\psi(x)dx + \frac{1}{s(s-1)}\;,
\end{align*}
where in $(*)$ the substitution $x=\tfrac{1}{y},dx = -\tfrac{1}{y^2}dy,\int_0^1\rightarrow \int_\infty^1$ was used.
Recall that we started with $\pi^{-\tfrac{s}{2}}\Gamma\left(\tfrac{s}{2}\right)\zeta(s) = \int_0^\infty x^{\tfrac{s}{2}-1}\psi(x) dx$ and arrived at
\begin{equation}
\pi^{-\frac{s}{2}}\Gamma\left(\frac{s}{2}\right)\zeta(s) = \int_1^\infty\left( x^{\frac{s}{2}} + x^{\frac{1-s}{2}} \right) \frac{\psi(x)}{x}dx - \frac{1}{s(s-1)} \, . \label{eq:ZetaPsiIntegralRelation}
\end{equation}
Note that the last term carries the pole of $\Gamma$ at $s=0$ and of $\zeta$ at $s=1$. Note further that the r.h.s. does not change under $s\rightarrow 1-s$, which implies
Riemann's functional equation
\begin{equation*}
\pi^{-\frac{s}{2}}\Gamma\left(\frac{s}{2}\right)\zeta(s) = \pi^{-\frac{1-s}{2}}\Gamma\left( \frac{1-s}{2}\right) \zeta(1-s)\; .
\end{equation*}
Riemann used 4-5 lines to derive this relation!
\end{proof}

In \eqref{eq:ZetaPsiIntegralRelation} we used
\begin{align*}
x^\frac{s}{2} = x^{\frac{\sigma + it}{2}} &= e^{\frac{\sigma\ln(x)}{2} + i\frac{t}{2}\ln(x)} = e^{\frac{\sigma\ln(x)}{2}}\left[\cos\left(\frac{t}{2}\ln(x)\right) + i\sin\left(\frac{t}{2}\ln(x)\right)\right]\\
x^{\frac{1-s}{2}} &= e^{\frac{(1-\sigma)\ln(x)}{2}}\left[\cos\left(\frac{t}{2}\ln(x)\right)-i\sin\left(\frac{t}{2}\ln(x)\right) \right]\\
x^\frac{s}{2}+x^{\frac{1}{2}(1-s)} &= \left(e^{\frac{\sigma \ln(x)}{2}} + e^{\frac{(1-\sigma)\ln(x)}{2}}\right) \cos\left(\frac{t}{2}\ln(x)\right)\\
&\overset{y=\frac{t}{2}\ln(x)}{=} \left(e^{\sigma\frac{y}{2}} + e^{(1-\sigma)\frac{y}{2}}\right)\cos(y)\\
&\overset{R.H.: \sigma = \frac{1}{2} }{=} \left( e^\frac{y}{2t} +e^\frac{y}{2t} \right)\cos(y) = 2 e^\frac{y}{2t} \cos(y)\\
&= 2e^{\frac{1}{4}\ln(x)}\cos(y) = 2x^\frac{1}{4}\cos\left(\frac{t}{2}\ln(x)\right)
\end{align*}
and whose imaginary part vanishes for $\sigma = \tfrac{1}{2}$. Thus
\begin{equation}
 \Xi(t) := \xi\left(\frac{1}{2}+it\right) = \frac{1}{2}+\frac{1}{2}s(s-1)\int_1^\infty \psi(x) \cdot 2 \cdot e^{\frac{1}{4}\ln(x)}\cos\left(\frac{t}{2}\ln(x)\right)\frac{dx}{x}
\end{equation}
is a real function, which is mentioned in Riemann's Berlin paper on p.147 as
\begin{equation}
\Xi(t) = \frac{1}{2}-(t^2+\frac{1}{4})\int_1^\infty \psi(x)x^{-\frac{3}{4}}\cos\left(\frac{t}{2}\ln(x)\right) dx \, ;
\end{equation}
furthermore,
\begin{equation}
\Im\xi\left(\frac{1}{2}+it\right) = 0,\;\Rightarrow\; \xi\left(\frac{1}{2}+it\right) = \Xi(t) \in\mathbb{R}\;.
\end{equation}
\section*{What is a function?}

Why is $1+2+3+4+\dots = -\tfrac{1}{12}$ a regularized value?
A normal reaction to this result:

This is not a true result. It is hogwash to say that $1+2+3+\dots $ has a finite value, as long as one does not specify what a function is 
(the concept of a function) and how it is calculated, i.e., which representation is chosen, what its domain of definition is, etc.

The following two statements are, however, true:
\begin{align*}
1+2+3+4+\dots &\rightarrow \infty,\quad \text{i.e., divergent}\\
\zeta_{\text{Riemann}}(-1) &= -\frac{1}{12}\quad .
\end{align*}

\textbf{Question}: In which representation is the latter statement true? We need a more general understanding of a function as well as the representation
in which the value of the function is calculated.

It is well known that a function can have several different representations, e.g., taking the sine function:
\begin{equation}
f(z) =\begin{cases} \sin (z) & \\
 \frac{e^{iz}-e^{-iz}}{2i} & \text{Euler} \\
 z -\frac{z^3}{3!}+\frac{z^5}{5!}-\dots & \text{Taylor expansion}\\
  z\prod_{n=1}^\infty \left(1-\frac{z^2}{\pi^2n^2}\right) & \text{product expansion}
  \end{cases}
\end{equation}

The Taylor expansion is an infinite-sum expansion of the sine function, one needs only powers of $z$.
The product expansion of the sine function needs all the infinitely many zeroes of the sine function.
One sees that there are \underline{many}  different ways to write a single function (e.g., sine), i.e., many different expressions
for performing various calculations !

What does all of this mean for the zeta function? Let's start with Euler's definition (1737):
\begin{align*}
\zeta(s) &= \sum_{n=1}^\infty \frac{1}{n^s},\quad s>1\\
&= 1+\frac{1}{2^s}+\frac{1}{3^s}+\dots,\quad s>1\;\text{for convergence}\ ,
\end{align*}
which is a sum of reciprocal powers of integers.
Evidently substituting \underline{negative} numbers for $s$ is not allowed, not even $s=1$ is permitted.

If one ignores the convergence condition $s>1$, then one can write
\begin{equation}
\zeta_{\text{Euler}}(-1) = 1+\frac{1}{2^{-1}}+\frac{1}{3^{-1}}+\dots = 1+2+3+4+\dots\quad ,
\end{equation}
which is pure nonsense, because it is not correctly defined. $s=-1$ is simply not allowed in Euler's definition (representation)
of the zeta function, which is only defined on the real axis $1<x\equiv s$. But there is another representation attributed to Riemann, which 
can be extended into the whole complex plane, $s\in\mathbb{C}\setminus\left\lbrace 0,1\right\rbrace$, i.e., including the value $\Re(s) = -1$.
\begin{equation}
\zeta(s) = \begin{cases}
\zeta_E(s) = \sum_{n=1}^\infty \frac{1}{n^s} & \Re(s)>1,\text{Euler (1797)}\\
\zeta_R(s) = 2^s\pi^{s-1}\sin\left(\frac{\pi s}{2}\right)\Gamma(1-s)\zeta(1-s) & s\in\mathbb{C}\setminus\lbrace 0,1\rbrace, \text{Riemann (1859)}
\end{cases}
\end{equation}
Note that the latter function is not given as a series but as a meromorphic function.

In Riemann's representation we obtain
\begin{align*}
\zeta_R(-1) &= 2^{-1}\pi^{-2}\sin\left(\frac{-\pi}{2}\right)\Gamma(1-(-1))\zeta(1-(-1))\\
&=2^{-1}\pi^{-2}(-1)\Gamma(2)\zeta(2)\\
&= 2^{-1}\pi^{-2}(-1)\cdot 1\cdot \frac{\pi^2}{6} = -\frac{1}{12}\quad,
\end{align*}
where in the third equality we used $\Gamma(2) = (2-1)!\cdot 1 = 1,\; \zeta(2) = 1+\tfrac{1}{2^2}+\tfrac{1}{3^2}+\dots = \tfrac{\pi^2}{6}$.

This is a true statement in Riemann's zeta-function representation
\begin{equation}
-\frac{1}{12} = \zeta_R(-1) \neq \zeta_E(-1) = \sum_{n=1}^\infty \frac{1}{n^s}\vert_{s=-1} \equiv 1+2+3+4+\dots
\end{equation}
whereas Euler's representation is not defined for $s=-1$.

The prime number counting function $\pi(x)$.

\textbf{Claim}:
\begin{align}
\frac{\ln\zeta(s)}{s} &= \int_2^\infty \frac{\pi(x)}{x(x^s-1)}dx,\quad s>1\label{eq:ZetaIntegralRepresentation}\\
\zeta(s) &= \prod_{p\in\text{primes}}\frac{1}{1-p^{-s}},\quad s>1\notag\\
\ln\zeta(s) &= \ln\prod_{p\in\text{primes}}\frac{1}{1-p^{-s}} = \sum_{p\in\text{primes}}\ln\frac{1}{1-p^{-s}}\notag
\end{align}
where $\pi(x)$ is the number of primes smaller than $x$. Replacing the summation over the primes by a summation over all integers yields
\begin{align}
\ln\zeta(s) &= \sum_{n=2}^\infty \lbrace \pi(n)-\pi(n-1)\rbrace\ln\frac{1}{1-n^{-s}}\quad \label{eq:LogZeta}
\intertext{where}
\pi(n) -\pi(n-1) &= \begin{cases}
1, & n\in\text{primes}\\
0, & \text{else}
\end{cases}\notag
\intertext{projects out the primes, e.g.,}
\pi(2)-\pi(1) &= 1-0 = 1\notag\\
\pi(3)-\pi(2) &= 2-1 = 1\notag\\
\pi(4)-\pi(3) &= 2-2 = 0\notag\\
\vdots &\notag\\
\eqref{eq:LogZeta} \Rightarrow \ln\zeta(s) &= \sum_{n=2}^\infty \pi(n) \ln\frac{1}{1-n^{-s}} - \sum_{n=2}^\infty\pi(n-1)\ln\frac{1}{1-n^{-s}}\notag\\
&= \sum_{n=2}^\infty \pi(n) \ln\frac{1}{1-n^{-s}} - \sum_{n=2}^\infty\pi(n)\ln\frac{1}{1-(n+1)^{-s}}\notag\\
&= \sum_{n=2}^\infty \pi(n)\left(\ln(1-(n+1)^{-s})-\ln(1-n^{-s})\right) . \label{eq:LogZetaExpanded}
\end{align}
Now use
\begin{equation}
\frac{d}{dx}\ln(1-x^{-s}) = \frac{1}{1-x^{-s}}(sx^{-s-1}) = \frac{s}{x(x^s-1)}.
\end{equation}
Integrate both sides to obtain 
\begin{equation}
\ln(1-x^{-s}) = s\int \frac{1}{x(x^s-1)}dx + C
\end{equation}
and use it in \eqref{eq:LogZetaExpanded}, whilst converting the indefinite integral into one over $[n,n+1]$:
\begin{align*}
\ln\zeta(s) &= \sum_{n=2}^\infty \underbrace{\pi(n)}_{\text{const. under integral}}\int_{n}^{n+1} \frac{s}{x(x^s-1)}dx \\
&= \sum_{n=2}^\infty \int_n^{n+1} \frac{s\pi(x)}{x(x^s-1)}dx\quad n: 2\rightarrow 3,3\rightarrow4,\dots\\
\ln\zeta(s) &= \int_2^\infty \frac{s\pi(x)}{x(x^s-1)}dx
\intertext{or}
\frac{\ln\zeta(s)}{s} &= \int_2^\infty \frac{\pi(x)}{x(x^s-1)}dx\quad .
\end{align*}
This concludes the proof.

For $s>1$ there are no non-trivial zeroes of $\zeta$.  Such are located in the critical strip $0<\Re(s)=\sigma<1$. The Riemann Hypothesis states that $\sigma=\tfrac{1}{2}$ for all zeroes of the $\zeta$ function.
Hence the formula \eqref{eq:ZetaIntegralRepresentation} is not applicable and we have to make an analytic continuation into the entire complex $s$ plane.

\end{document}